\documentclass[preprint,12pt]{elsarticle}
\usepackage[utf8]{inputenc}
\usepackage[T1]{fontenc}
\usepackage{hyperref}
\usepackage{amsmath}
\usepackage{listings}
\usepackage{todonotes}
\usepackage{multirow}
\usepackage{float}
\usepackage{booktabs}

\providecommand{\abs}[1]{\lvert#1\rvert}

\newcommand{\powipmi}{\pi_{1+}^{-1}}
\newcommand{\expmi}{\exp_{-1}}
\newcommand{\logip}{\log_{1+}}

\let\code\texttt

\lstnewenvironment{lstfloat}[1][]%
{\lstset{basicstyle=\small\ttfamily,basewidth=0.5em,frame=TBlr,float=tbp,#1}}{}
\lstnewenvironment{lstdisplay}[1][]%
{\lstset{basicstyle=\ttfamily,xleftmargin=1em,basewidth=0.5em,frame=none,#1}}{}

\journal{Computers and Fluids}

\begin{document}

\begin{frontmatter}

\title{Computational performance versus numerical accuracy in the acoustic terms of weakly-compressible SPH}

\author[ingv]{Giuseppe Bilotta\corref{cor}}
\ead{giuseppe.bilotta@ingv.it}
\cortext[cor]{corresponding author}

\author[uw]{Alireza Zarei}
\author[ingv]{Roberta Cristofaro}
\author[cnam]{Alexis Hérault}
\author[cea]{Élie Saikali}
\author[uw]{Morteza Derakhti}

\affiliation[ingv]{
	organization={Osservatorio Etneo, Istituto Nazionale di Geofisica e Vulcanologia},
	city={Catania},
	country={Italy}
}

\affiliation[uw]{
	organization={University of Washington},
	city={Seattle},
	state={WA},
	country={USA}
}

\affiliation[cnam]{
	organization={Laboratoire Mod\'elisation math\'ematique et num\'erique, Conservatoire National des Arts et M\'etiers},
	city={Paris},
	country={France}
}

\affiliation[cea]{
	organization={Université Paris-Saclay, CEA, Service de Génie Logiciel pour la Simulation, 91191, Gif-sur-Yvette},
	country={France}
}

\begin{abstract}
Accurate computation of the acoustic terms in weakly-compressible SPH, including the equation of state, its inverse, and the evaluation of the local speed of sound, is essential to ensure energy conservation, robust solid boundary conditions, and correct open boundary behavior.
In the classic Cole's equation of state, using values of the polytropic constant $\gamma \ne 1$ can lead to significant precision loss,
particularly when computations are performed in single precision, as is common in GPU implementations.
We extend here our previous results on the optimal choice of state variable by introducing a new higher-precision form for all three acoustic terms,
denoted \texttt{explog}, and assess its numerical accuracy and computational performance across a variety of devices and test cases.
We find that the new form evaluates the acoustic terms with a maximum ULP of~2 for arguments less than $0.125$,
and has higher errors in less than 1\% of the arguments larger than the threshold.
Combined with a rewritten expression of the pressure at solid boundary particles in the ``dummy'' model,
it provides consistently more accurate results, maintaining deviations from the reference pressure profile to below 2\%,
even in presence of a rough bottom terrain.
Computationally, the new formulation is also faster than the naive one on GPU, and only marginally slower on CPU.
\end{abstract}

\begin{highlights}
\item Weakly-compressible SPH is sensitive to noise in the density and pressure fields;
\item The standard form of the Cole's equation of state is not numerically robust for small fluctuations;
\item Reformulating the acoustic terms in terms of \code{expm1} and \code{log1p} improves numerical robustness;
\item The new \code{explog} form has a ULP $\le 2$ for arguments less than 1/8;
\item The new formulation is faster than less accurate choices on GPU, but not on CPU.
\end{highlights}


\end{frontmatter}

\section{Introduction}

The weakly compressible form of Smoothed Particle Hydrodynamics (WCSPH) is a meshless Lagrangian method
well suited for simulating violent flows~\cite{monaghan1994simulating},
owing to its natural and implicit treatment of the free surface, including splashing, fragmentation, and merging \cite{ monaghan2005smoothed,gomez2010state}.

Coupled with an explicit integration scheme, WCSPH is also well suited for implementation on high-performance parallel computing hardware,
owing to the natural mapping between numerical nodes (particles) and computational work-items. This makes it particularly compatible with the stream processing paradigm commonly used to program Graphic Processing Units (GPUs).
Indeed, even before the emergence of general-purpose GPU computing (GPGPU), researchers had already begun to exploit the computational power of these accelerators
to model SPH and related particle systems \cite{takahashi2003real,kipfer2004uberflow}. 

Although WCSPH as a numerical method and GPUs as computational hardware are well suited for each other,
several aspects must be carefully evaluated, as a naive implementation can lead to numerical instabilities or poor performance%
~\cite{herault_accuracy_2014,bilotta_design_2018,saikali_2020,sph_gpu_trap,fast_neib_list}.

One of the key challenges lies in numerical precision.
Since GPUs were originally designed for graphics applications,
most consumer-grade devices offer very high single-precision but very low double-precision
floating-point throughput.
The performance ratio between single and double precision typically ranges from 1:8 to 1:32 on consumer GPUs, which can significantly reduce the computational advantages of GPU acceleration.
Even many high-end workstation or server GPUs do not reach the 1:2 single-to-double precision performance ratio typical of modern CPUs and are instead limited to ratios of about 1:3 or 1:4,
despite their considerably higher cost compared to consumer GPUs.

Implementations of WCSPH for GPUs must therefore balance performance and numerical accuracy. 
They can either limit computations to single precision, with a potential loss of stability,
or adopt double precision for higher accuracy, possibly restricting its use to the most sensitive parts of the code.
The latter option, however, often results in a significant drop in performance on lower-end consumer GPUs~\cite{dominguez-neibs-list,dominguez_double_precision}.

A practical compromise is to rely on single-precision computations only and apply numerical tricks that maximize the effective use of the available precision.
Such techniques help mitigate rounding and cancellation errors that would otherwise accumulate in single-precision arithmetic. 
They range from classic approaches, such as compensated summation~\cite{kahan1965,neumaier1974,klein2006}, to appropriate choices of formulations and reference variables~\cite{herault_accuracy_2014,saikali_2020},
and may even extend to redesigning well-established algorithms to avoid catastrophic cancellations~\cite{bilotta_bicgstab}.

Ultimately, all these approaches share the common feature of requiring some degree of transformation of the analytical formulation,
thereby breaking the immediate correspondence between analytical expressions and their implementation.
This contrasts with formulations that are designed, already at the analytical level, to reduce the numerical noise that develops in numerical discretizations regardless of precision.
Examples include artificial viscosity~\cite{monaghan_2005},
density diffusion~\cite{monaghan1983shock,molteni2009density,molteni_colagrossi_2009,antuono_2011,marrone_2011,green_2019},
and kernel gradient corrections~\cite{vila_1999,chen_1999},
which are generally introduced a posteriori in the (analytical) discretized equations,
sometimes sacrificing strict consistency in favor of improved convergence and stability~\cite{chen_beraun_carney_corrective,bonet_variational,zago_ccsph},
as discussed recently also in~\cite{letouze_review_2025}.

In WCSPH, the acoustic terms are particularly important at both the computational and analytical levels.
Because pressure is computed explicitly from density through an equation of state,
numerical fluctuations in the density field can propagate to the velocity field and accumulate over time, undermining the stability of the simulation. 
This makes it extremely difficult to reproduce simple hydrostatic cases (e.g., still water in a parallelepipedal box).
Although these effects are generally of second order in highly dynamic problems,
they can still have a negative impact on solid boundary models~\cite{adami_2012,mayrhofer_2015}
and on the implementation of open boundary conditions~\cite{ferrand_2017,tafuni_2018}.

This behavior is generally independent of the numerical precision adopted in the code,
although single-precision computations may lead to such instabilities after fewer iterations.
The importance of maintaining a smooth density (and pressure) field has led to the development of several density diffusion schemes in WCSPH,
ranging from periodic filtering or reconstruction methods~\cite{colagrossi_landrini_2003,seo2021}
to diffusion terms directly incorporated into the mass conservation equation~\cite{molteni2009density,molteni_colagrossi_2009,marrone_2011,deltaplus,green_2019,mr_deltaplus},
drawing parallels with the von Neumann artificial viscosity used in the equation of motion~\cite{monaghan_2005}.
A similar approach can also be used to include damping terms in the viscous contribution,
specifically designed to limit the formation of spurious acoustic waves~\cite{acoustic_damper}
which have been shown to be a significant issue in the application of WCSPH to engineering applications~\cite{meringolo_2017}.
These issues are further compounded by the contributions to the pressure field coming from the terms needed to model solid boundary conditions,
a challenging problem on its own merit within SPH~\cite{grand_challenges_2020}.

The work we present here is an analysis of the acoustic terms in WCSPH associated with the classic Cole's equation of state traditionally used in SPH.
We examine both the numerical accuracy and computational performance of the equation in its multiple implementable forms
and use this analysis to identify possible numerical optimizations for evaluating the pressure term in common WCSPH formulations, with or without density diffusion, depending on the chosen boundary model.

The paper briefly introduces the SPH method with a focus on the acoustic terms (Section~\ref{sec:sph_intro}).
This is followed by a study of the abstract analytical and computational aspects of the evaluation of the acoustic terms (Sections~\ref{sec:numerical_considerations} and~\ref{sec:computational_considerations}),
discussed independently of their specific application to SPH, since these results are of general interest to any numerical code in which the same functional forms appear.
Section~\ref{sec:dummy} then focuses on the application to SPH and presents a new functional representation for the ``dummy'' boundary model, leveraging the improved accuracy to achieve more stable and consistent numerical results.

\section{Weakly-compressible SPH}\label{sec:sph_intro}

We present here a quick introduction to SPH and its formalism. Further details can be found e.g. in~\cite{monaghan_2005,violeau_book}.

In SPH, the continuum is discretized as a set of particles, each carrying the physical properties (mass, density, volume, temperature, etc)
of a virtual volume of fluid.
Scalar and vector fields are approximated via convolution with a family of radial \emph{smoothing kernel}s $W(\cdot, h)$
that approximate Dirac's $\delta$ distribution, converging to it in the limit as the \emph{smoothing parameter} $h$ goes to~0,
and the integral of the convolution is then approximated with a summation over particles:
\[
f(\vec x) = \sum_j f(\vec x_j) W(\vec x - \vec x_j, h) V_j
\]
where $V_j$ is the volume of particle $j$.
Using smoothing kernels with compact support allows such summations to be limited to neighboring particles only.

Gradients can be approximated leveraging Green's theorem to ``shift'' the graident operator from the field to the smoothing kernel:
\[
\nabla_x f(\vec x) = \sum_j f(\vec x_j) \nabla_x W(\vec x - \vec x_j, h) V_j + B(\vec x)
\]
where $B(\vec x)$ is a boundary term. When the smoothing kernel has compact support, we have $B(\vec x) = 0$ everywhere
except for points $\vec x$ whose distance to the boundary is less than the kernel support radius.
The kernel gradient expression can be symmetrized in a number of way, e.g. by adding or subtracting $0 = f(\vec x) \nabla 1$,
resulting in a variety of expressions for gradients in SPH.

Due to kernel symmetry, we can write $\nabla_x W(\vec x - \vec x_j, h) = (\vec x - \vec x_j) F(\vec x - \vec x_j, h)$
where $F(r, h) = (1/r) \partial W(r, h)/\partial r$.
This form improves numerical stability when $F$ has an analytic representation where the $1/r$ factor cancels out with radial terms
in the derivative of~$W$.

Using the common SPH shorthand conventions $\vec x_{ij} = \vec x_i - \vec x_j$,
$W_{ij} = W(\vec x_{ij}, h)$ and $F_{ij} = F(\vec x_{ij}, h)$,
and $f_{ij} = f(\vec x_i) - f(\vec x_j)$ for any field $f$,
the standard mass and momentum continuity equations in a weak compressibility regime can be discretized in SPH
following Monaghan's ``golden rules''~\cite{monaghan_2005}:
\begin{align}
\frac{D\rho_i}{Dt} &= -\sum_j m_j \vec u_{ij} \cdot \vec x_{ij} F_{ij}\label{eq:mass-cont-sph1}\\
\frac{D \vec u_i}{Dt} &= -\sum_j \left(
\left( \frac{P_i}{\rho_i^2} + \frac{P_j}{\rho_j^2}\right) \vec x_{ij}
+ 2K \frac{m_j}{\rho_i \rho_j} \bar\mu_{ij}
\frac{\vec x_{ij} \cdot \vec u_{ij}}{\vec x_{ij}^2} \vec x_{ij} \right) F_{ij} + \vec g,\label{eq:ns-sph1}
\end{align}
optionally with additional terms for artificial viscosity and\slash or density diffusion;
in the above equations, $\vec u_i$ represents the velocity of particle $i$,
$m_i, \rho_i, P_i$ its mass, density and pressure,
$\bar \mu_{ij}$ the average viscosity between particles $i, j$,
and $K$ is a constant that depends on the problem dimension and choice of viscosity averaging function~\cite{monaghan_2005,morris_1997,violeau_book}.

\subsection{Acoustic terms}

The weak compressible formulation of SPH closes the momentum equation with a functional relation between pressure and density
under a barotropic assumption.
Commonly used is Cole's equation of state~\cite{cole,batchelor}
(also known as the Tait--Murnaghan equation or somewhat improperly simply the Tait equation~\cite{macdonald_1966}):
\begin{equation}\label{eq:eos}
P(\rho) = B \left( \left(\frac{\rho}{\rho_0}\right)^\gamma -1\right)
\end{equation}
where $\rho_0$ is the at-rest density, $B$ the bulk modulus and $\gamma$ the polytropic constant.
Conversely, given the pressure $P$, the corresponding density can be derived as
\begin{equation}\label{eq:inverse_eos}
\rho(P) = \rho_0 \left(1 + \frac{P}{B}\right)^{1/\gamma}.
\end{equation}

The speed of sound can be obtained from $c(\rho) = \sqrt{\partial P/\partial \rho}$:
\begin{equation}\label{eq:soundspeed}
c(\rho) = \sqrt{\frac{B \gamma}{\rho_0}} \left(\frac{\rho}{\rho_0}\right)^{\frac{\gamma-1}{2}}
\end{equation}
and letting $c_0$ be the sound speed for the fluid at rest, we get
\begin{equation}\label{eq:bulk}
B = \frac{\rho_0 c_0^2}{\gamma}.
\end{equation}

Since WCSPH is frequently used with explicit integration schemes,
for which the speed of sound is often the dominant factor
(unless simulating highly viscous fluids at very fine resolutions~\cite{bilotta_bicgstab}),
the physical speed of sound is rarely used,
preferring the adoption of a fictitious sound speed.
To ensure the weakly compressible regime is valid, density variations should be
less than $1\%$, for which it is generally sufficient to take $c_0$ at least 10 times higher than the maximum flow velocity,
including the so-called hydrostatic limit $\sqrt{gH}$, with $H$ the reference height,
in order to avoid the collapse of a column of water~\cite{violeau_book}.

\subsection{Boundary conditions}\label{sec:dummy_intro}

Near the boundary of the fluid domain, the SPH interpolations are defective due to missing contributions to the kernel support
from outside the domain.
Near the free surface, these missing contributions can generally be ignored assuming that the field values are null,
but near solid boundaries, corrective terms must be applied.
Among the many forms that boundary conditions can taken in SPH~\cite{Liu:2003,Monaghan:2009,violeau_book},
the most common ones are based on a particle discretization similar to the fluid one,
with several layers of boundary particles included in the simulation
in order to complete the kernel support for fluid particles near the boundary.

The simplest approach (so-called dynamic boundary conditions~\cite{dalrymple2001sph,Crespo:2007}) has boundary particles
behave like fluid particles, evolving the density according to the mass continuity equation,
but prescribing velocity based on the rigid body motion (if any) of the boundary.
While this is generally sufficient to avoid fluid particles penetrating the boundary,
it does not provide the correct stress and strain at the boundary, and can produce artifacts
such as a `lift' of the head in a bore, or particles sticking to the wall when fluid retreats.%

A more sophisticated approach is the ``dummy'' boundary model introduced in~\cite{adami_2012}.
In this case, boundary particle properties are computed in order to prescribe both the normal stress
and the no-slip conditions.
In particular, the boundary particle pressure is obtained from an interpolation of the adjacent fluid pressure,
with a correction to compensate for the hydrostatic term:
\begin{equation}
P_w = \frac{
\sum_j (P_j + \rho_j \vec g \cdot \vec x_{wj}) W(\vec x_{wj}, h)
}{
\sum_j W(\vec x_{wj}, h)
}.\label{eq:dummy_pressure}
\end{equation}
The particle density can then be computed from the pressure using the inverse of the equation of state,
as needed.

It should be apparent that the accuracy in the computation of the wall pressure $P_w$
has a significant impact on the accuracy of this boundary condition term,
since incorrect evaluations may lead to the wall attracting (if the pressure is underestimated)
or repulsing (if the pressure is overestimated) fluid particles,
leading respectively to possible penetration, or to an underestimation of the wall friction term.
Improving the evaluation of these terms is the focus of this paper.

\section{Numerical considerations on the acoustic term evaluation}\label{sec:numerical_considerations}

When working with single-precision (32-bit) floating point, with machine epsilon
\[
\epsilon_{M,32} = 2^{-23} \approx 1.19209289551\cdot10^{-7},
\]
manipulating quantities such as density, pressure and the sound speed of the fluid can lead to low accuracy in the results,
due to the amount of precision consumed by the magnitude of the physical quantities.
Remember that the machine epsilon is defined as the value of the least significant bit (ULP, Unit of Last Position or equivalently
Unit of Least Precision) in the floating-point representation of the value $1$.
This effectively defines the ratio between the least and most significant digit in any representable number,
and the smallest relative change in a value that can be correctly represented.

Now consider e.g.~$\rho_0$ in the order of $10^3$ (in SI units of $\text{kg}\,\text{m}^{-3}$) for liquids:
this effectively consumes 3 decimal digits of precision,
so that only variations in the order of $\epsilon_M \rho_0$ are representable, which in this case means
variations in the order of $10^{-7}\cdot10^3 = 10^{-4}$ in single precision.
If we then take $c_0$ in the order of $10$ (in SI units of $\text{m}\,\text{s}^{-1}$),
$c_0^2$ will consume 2 decimal digits, and $B$ will be in the order of $10^6$,
leaving only one decimal digit for representable absolute variations in the pressure
($O(\delta P) = O(\epsilon_M P) = 10^{-7}\cdot10^6 = 10^{-1}$).

This loss of precision can lead to numerically unstable simulations even for analytically robust formulations,
leading to effects ranging from spurious explosions to fluid particles unphysically penetrating the solid boundary layers~\cite{rustico_advances_2014,herault_accuracy_2014}.
Although adimensionalization can help by bringing function arguments within the ranges of higher precision,
it does not necessarily improve the representability of variations, since we're still within
machine epsilon of~1.

Among the several numerical strategies that can be adopted to compensate for the loss of precision,
the most important for the treatment of the acoustic terms is the shift from the physical density $\rho$ to the \emph{relative density variation}
\cite{herault_accuracy_2014,saikali_2020}
(also known as the \emph{condensation} in acoustic wave theory~\cite{fundamental_acoustics}),
defined as:
\[
\tilde\rho = \frac{\rho}{\rho_0} - 1.
\]

In a weakly-compressible regime, $\abs{\tilde\rho} \ll 1$, and this choice allows tracking of density fluctuations as small as the smallest normal number
\code{FLT\_MIN}
\begin{equation}
f_{m,32} = 2^{-126} \approx 1.17549435082228751 \cdot 10^{-38} \label{flt_min}
\end{equation}
without loss of precision.
It is thus recommended that $\tilde\rho$ be chosen as state variable, with $\rho = \rho_0(\tilde\rho+1)$ a derived quantity.
The mass continuity equation for example becomes:
\begin{equation}\label{eq:mass-cont-sph-tilde}
\frac{D\tilde\rho_i}{Dt} = -\sum_j \frac{m_j}{\rho_{0,i}} \vec u_{ij} \cdot \vec x_{ij} F_{ij}
\end{equation}
and the acoustic terms can be rewritten as
\begin{equation}\label{eq:eos_rhotilde}
P(\tilde\rho) = B \left( (1 + \tilde\rho)^\gamma - 1\right)
\end{equation}
for the equation of state, with inverse
\begin{equation}\label{eq:inverse_eos_rhotilde}
\tilde\rho(P) = \left( 1 + \frac{P}{B}\right)^{1/\gamma} - 1,
\end{equation}
and
\begin{equation}\label{eq:c_rhotilde}
c(\tilde\rho) = c_0 (1 + \tilde\rho)^{(\gamma - 1)/2}.
\end{equation}
for the sound speed.

If we introduce the scaled density $\tilde P = P/B$ the equation of state and its inverse show a beautiful symmetry:
\begin{align}
\tilde P(\tilde\rho) & {}= (1 + \tilde\rho)^{\gamma\phantom{/1}} - 1,\label{eq:ptilde_rhotilde}\\
\tilde \rho(\tilde P) & {}= (1 + \tilde P)^{1/\gamma} - 1\label{eq:rhotilde_ptilde}
\end{align}
that, as we will show in Section~\ref{sec:dummy}, will allow us to improve the numerical accuracy and stability of simulations.
We remark that under the condition $-1 < \tilde\rho < 2^{1/\gamma} - 1$ we also have $\abs{\tilde P} < 1$,
since
\begin{align*}
\abs{ \tilde P } < 1 &\iff\\
\abs{ (1+\tilde\rho)^\gamma - 1} < 1 &\iff\\
-1 < (1+\tilde\rho)^\gamma - 1 < 1 &\iff\\
0 < (1+\tilde\rho)^\gamma < 2 &\iff\\
0 < 1+\tilde\rho < 2^{1/\gamma} &\iff\\
-1 < \tilde\rho < 2^{1/\gamma} - 1.  &
\end{align*}

\section{Computational considerations}\label{sec:computational_considerations}

Equations \eqref{eq:ptilde_rhotilde}, \eqref{eq:rhotilde_ptilde} and \eqref{eq:c_rhotilde}
can all be expressed in terms of the function
\[
\powipmi(x, \alpha) = (1+x)^\alpha - 1,
\]
since we can write:
\begin{align}
\tilde P(\tilde\rho)  & {}= \powipmi(\tilde\rho, \gamma),\label{eq:prhopow}\\
\tilde \rho(\tilde P) & {}= \powipmi(\tilde P, 1/\gamma),\label{eq:rhoppow}\\
c(\tilde\rho) & {}= c_0 \left( \powipmi(\tilde\rho, (\gamma - 1)/2) + 1\right).\label{eq:crhopow}
\end{align}
It may seem that expression \eqref{eq:crhopow} is more complex than the original form \eqref{eq:c_rhotilde},
but, as we shall see momentarily,
the adoption of $\powipmi$ leads to significant numerical improvements.

In code, we refer to this function as \code{pow1pm1},
inspired by the functions \code{expm1}, \code{log1p} defined by the IEEE-754 standard~\cite{ieee754}
and corresponding respectively to $\expmi(x) = e^x - 1$ and $\logip(x) = \log(1+x)$,
designed to assist developers in computing exponentials and logarithms with higher accuracy as $x\to0$.
In~fact, we have:
\[
\powipmi(x, \alpha) = (1+x)^\alpha - 1 = e^{\alpha \log(1+x)} - 1 =
\expmi(\alpha \logip(x)).
\]

The Taylor series expansion of $\powipmi(x, \alpha)$ around $x = 0$ is
\[
\powipmi(x, \alpha) = \alpha x \left(
1 + \sum_{k=1}^\infty \frac { \prod_{i=1}^k (\alpha - i) }{(k+1)!} x^k
\right)
\]
that can be written in Horner form~\cite{horner1819}
as
\begin{equation}
\powipmi(x, \alpha) = \alpha x \left(
1 + \frac{\alpha - 1}{2} x \left( 1 + \frac{\alpha - 2}{3} x \left(
1 + \cdots\phantom{\frac12} \right) \right) \right).
\label{eq:pow1pm1_horner}
\end{equation}

This fact was leveraged in~\cite{herault_accuracy_2014}
to optimize the computation of $P(\rho)$ and $P(\rho)/\rho^2$ in the case $\gamma = 7$,
where both expressions can be computed quickly and accurately
as polynomials in Horner form using fused multiply-add (FMA) operations.

\subsection{Computational speed versus accuracy}\label{sec:powforms}

The polynomial trick considered by~\cite{herault_accuracy_2014} can only be applied
when $\powipmi(x, \alpha)$ is a polynomial.
In the case of the three acoustic terms we are considering
\eqref{eq:prhopow}--\eqref{eq:crhopow}
(pressure from density, density from pressure, speed of sound)
with exponents respectively $\gamma, 1/\gamma, (\gamma - 1)/2$,
the only case in which this applies to all of them is for $\gamma = 1$.

The simplification in the expression of the acoustic terms for $\gamma = 1$
make it a common choice in SPH,
even though the exponent is considerably smaller than the physical polytropic constant of water,
$\gamma_w \approx 7.5$.
Outside of the case $\gamma = 1$, however (e.g. with the other common choice $\gamma = 7$),
not all the acoustic terms can have integer exponents at the same time,
and the polynomial simplifications of~\cite{herault_accuracy_2014} are only possible for some of them at most.
The issue becomes then: what is the fastest and most accurate way to evaluate $\powipmi$ in the non-polynomial
(or in the general) case?

Coding-wise, we can look at $\powipmi$ as having three possible functional forms,
illustrated here for the single-precision \code{float} type:
\begin{lstdisplay}
float pow1pm1(float x, float a) /* power form */
	{ return powf(1 + x, a) - 1; }
float pow1pm1(float x, float a) /* explog form */
	{ return expm1f(a*log1pf(x)); }
float pow1pm1(float x, float a) /* Taylor form */
	{ return /* omissis: Horner form Taylor expansion */ }
\end{lstdisplay}
Additionally, some programming languages and hardware combinations support
intrinsics such as \code{\_\_powf} that sacrifice accuracy for speed,
but can still be used as an alternative in the power form.

While these forms are analytically equivalent, computationally they differ significantly.
The power form, which is what is naively used in many codes, will return the value 0 if $\abs{x} < \epsilon_M$,
since for such small values of $x$ we have $1 \oplus x = 1$ (indicating by $\oplus$ the correctly-rounded result
of the operation between two floating-point values).

Note that increasing the precision will simply change the threshold below which the actual value of $x$ will be ignored:
$\epsilon_{M,32}  = 2^{-23}$ for single precision,
$\epsilon_{M,64}  = 2^{-52}$ for double,
$\epsilon_{M,80}  = 2^{-63}$ in the extended 80-bit format supported as `long double' on some architectures,
$\epsilon_{M,128} = 2^{-112}$ for the 128-bit quadruple precision.
Note that
even \emph{with computations done in quadruple precision}
the power form will \emph{fail to return the correct value for single-precision} \code{FLT\_MIN} \eqref{flt_min}.

The explog form, on the other hand, can be expected to produce very accurate results,
assuming the implementations of \code{expm1f} and \code{log1pf} are correctly rounded,
and thus produce results within 1 ULP of the exact value.
In this case, the argument to $\expmi$ will be accurate to within 1 ULPs,
and given the magnitude of the arguments, this will result in an overall error of at most 2 ULPs
in the final result.
This high accuracy allows us to use the explog in e.g. double or extended precision as reference value
for the single precision function in our validation tests.

\begin{lstfloat}[caption={Recursive function template definition for the Horner form
of the Taylor expansion of the $\powipmi(x, \alpha)$ function.
The \code{m*x*prev + 1} can be decomposed in three different ways as argument to the
fused multiply-add: in addition to the one illustrated in code, the forms
\code{fmaf(m*prev, x, 1.0f)} and \code{fmaf(x*prev, m, 1.0f)} are also possible.},
label={lst:pow1pm1_horner}]
template<unsigned order>
float taylor_form(float x, float a, float prev=1.0f);

template<>
float taylor_form<1>(float x, float a, float prev)
{ return a*x*prev; }

template<unsigned order>
float taylor_form(float x, float a, float prev)
{
	float m = (a - (order - 1))/order;
	return taylor_form<order-1>(x, a, fmaf(m*x, prev, 1.0f));
}
\end{lstfloat}

The accuracy of the Taylor expansion depends on the magnitude of the arguments
and the order chosen for the expansion.
The Horner form~\eqref{eq:pow1pm1_horner} can be mapped easily to a recursive function template in C++
(Listing~\ref{lst:pow1pm1_horner}), allowing the selection of the Taylor expansion order
based on the preferred compromise between performance and accuracy.
Note that the Horner form contains \emph{two} multiplications and an addition for each nesting;
mapping them to FMA therefore requires a choice about
which multiplication should be done before the FMA and which one should be done inside,
leading to at least three different forms for the Taylor expansion:
for them we can expect minor differences in their results, as well as some difference in computational
performance due to a different degree of instruction-level parallelism.
We denote them by ``mx'' for \code{fmaf(m*x, prev, 1.0f)},
``mprev'' for \code{fmaf(m*prev, x, 1.0f)} and
``xprev'' for the \code{fmaf(x*prev, m, 1.0f)} form.

\subsection{Experimental results}

We evaluate both the numerical and computational performance of the different forms of $\powipmi(x, \alpha)$
with $x \in [0, 1[$ and a selection of the common values needed for $\alpha$,
corresponding to $\gamma, 1/\gamma$ and $(\gamma - 1)/2$ for $\gamma = 1, 3, 5, 7$.
The case $\alpha = 0$ is not presented here because all forms give the correct result in that case.
We also present the results for $\alpha = -1, -2$ as they are relevant to the pressure term computation
as explained later in Section~\ref{sec:dummy-pressure}.

The numerical results are evaluated by comparing the results from the single-precision forms
against the results obtained by the \code{long double} explog form.
The metrics we consider are the number of discrepancies
(values for which the single-precision results were different from the long double result),
the number of values for which the discrepancy is larger than 1 ULP (in single precision),
and the largest relative error in ULPs.

To account for the difference in error distribution for different magnitudes of $x$,
the $[0, 1[$ range has been split into three subranges: denormals $[0, f_{m,32}[ = [0, 2^{-126}[$,
``nano'' $[f_{m,32}, \epsilon_{M,32}[ = [2^{-126}, 2^{-23}[$,
``micro'' $[\epsilon_{M,32}, 2^{-6}[ = [2^{-23}, 2^{-6}[$,
``percent'' $[2^{-6}, 2^{-3}[$,
and ``perten'' $[2^{-3}, 1[$.
The marks at $2^{-6}=1.5625\cdot10^{-2}$ and $2^{-3} = 0.125$ were chosen
as the smallest powers of~2 that are not less than $1\%$ and $10\%$ respectively.

The code for these tests is available as auxiliary material to the manuscript.
The performance results refer to a machine with an AMD Ryzen 7 3700X processor
that can boost to 4.4GHz, running Linux 5.18 with vulnerability mitigations enabled.
The code was compiled with Clang 14, \code{-O3} and
machine-specific optimizations enabled (\code{-march=native}).
Due to the high optimization level, the compiler produces vectorized code
and interleaves instructions to maximize throughput where possible,
with additional performance benefits for some forms.
Additional performance comparisons in more practical contexts, running on GPUs
and including the faster and less accurate \code{\_\_powf} form will be presented
as part of the application to the dummy boundary model in Section~\ref{sec:dummy}.

\begin{table}
\caption{Performance, in $10^9$ elements (GE) per seconds, and accuracy of the $\powipmi$ function in the power, explog and Taylor forms 
in the $[0, 2^{-126}[$ range, covering $8,388,608$ single-precision floating-point values.}
\label{tab:pow1pm1.log.1}
\scriptsize\centering
\begin{tabular}{l|l lrrl}
$\alpha$ & method & GE/sec & discrepancies & ${}> 1$ ULP & max ULPs\\
\hline
\multirow{5}*{$-2$}
& power & 0.14 & 100\% & 100\% & ${} > 2^{23}$\\
& explog & 0.2 & 0\% & 0\% & 0\\
& Taylor mx & 2.8 & 0\% & 0\% & 0\\
& Taylor mprev & 2.1 & 0\% & 0\% & 0\\
& Taylor xprev & 2.1 & 0\% & 0\% & 0\\
\hline
\multirow{5}*{$-1$}
& power & 0.16 & 100\% & 100\% & ${} > 2^{22}$\\
& explog & 0.22 & 0\% & 0\% & 0\\
& Taylor mx & 2.8 & 0\% & 0\% & 0\\
& Taylor mprev & 2.1 & 0\% & 0\% & 0\\
& Taylor xprev & 2.1 & 0\% & 0\% & 0\\
\hline
\multirow{5}*{$1/7$}
& power & 0.16 & 100\% & 100\% & ${} > 2^{20}$\\
& explog & 0.22 & 0\% & 0\% & 0\\
& Taylor mx & 2.8 & 0\% & 0\% & 0\\
& Taylor mprev & 2.1 & 0\% & 0\% & 0\\
& Taylor xprev & 2.1 & 0\% & 0\% & 0\\
\hline
\multirow{5}*{$1/5$}
& power & 0.16 & 100\% & 100\% & ${} > 2^{20}$\\
& explog & 0.22 & 0\% & 0\% & 0\\
& Taylor mx & 2.8 & 0\% & 0\% & 0\\
& Taylor mprev & 2.1 & 0\% & 0\% & 0\\
& Taylor xprev & 2.1 & 0\% & 0\% & 0\\
\hline
\multirow{5}*{$1/3$}
& power & 0.16 & 100\% & 100\% & ${} > 2^{21}$\\
& explog & 0.22 & 0\% & 0\% & 0\\
& Taylor mx & 2.8 & 0\% & 0\% & 0\\
& Taylor mprev & 2.1 & 0\% & 0\% & 0\\
& Taylor xprev & 2.1 & 0\% & 0\% & 0\\
\hline
\multirow{5}*{$1$}
& power & 0.16 & 100\% & 100\% & ${} > 2^{22}$\\
& explog & 0.22 & 0\% & 0\% & 0\\
& Taylor mx & 2.8 & 0\% & 0\% & 0\\
& Taylor mprev & 2.1 & 0\% & 0\% & 0\\
& Taylor xprev & 2.1 & 0\% & 0\% & 0\\
\hline
\multirow{5}*{$2$}
& power & 0.16 & 100\% & 100\% & ${} > 2^{23}$\\
& explog & 0.21 & 0\% & 0\% & 0\\
& Taylor mx & 2.8 & 0\% & 0\% & 0\\
& Taylor mprev & 2.1 & 0\% & 0\% & 0\\
& Taylor xprev & 2.1 & 0\% & 0\% & 0\\
\hline
\multirow{5}*{$3$}
& power & 0.16 & 100\% & 100\% & ${} > 2^{23}$\\
& explog & 0.2 & 0\% & 0\% & 0\\
& Taylor mx & 2.8 & 0\% & 0\% & 0\\
& Taylor mprev & 2.1 & 0\% & 0\% & 0\\
& Taylor xprev & 2.1 & 0\% & 0\% & 0\\
\hline
\multirow{5}*{$5$}
& power & 0.16 & 100\% & 100\% & ${} > 2^{23}$\\
& explog & 0.2 & 0\% & 0\% & 0\\
& Taylor mx & 2.8 & 0\% & 0\% & 0\\
& Taylor mprev & 2.1 & 0\% & 0\% & 0\\
& Taylor xprev & 2.1 & 0\% & 0\% & 0\\
\hline
\multirow{5}*{$7$}
& power & 0.16 & 100\% & 100\% & ${} > 2^{23}$\\
& explog & 0.21 & 0\% & 0\% & 0\\
& Taylor mx & 2.1 & 0\% & 0\% & 0\\
& Taylor mprev & 1.7 & 0\% & 0\% & 0\\
& Taylor xprev & 2.1 & 0\% & 0\% & 0\\
\end{tabular}
\end{table}

\begin{table}
\caption{Performance, in $10^9$ elements (GE) per seconds, and accuracy of the $\powipmi$ function in the power, explog and Taylor forms
in the $[2^{-126}, 2^{-23}[$ range, covering $872,415,232$ single-precision floating-point values.}
\label{tab:pow1pm1.log.2}
\scriptsize\centering
\begin{tabular}{l|l lrrl}
$\alpha$ & method & GE/sec & discrepancies & ${}> 1$ ULP & max ULPs\\
\hline
\multirow{5}*{$-2$}
& power & 0.15 & 100\% & 100\% & $2^{24}$\\
& explog & 0.18 & 0.693\% & 0\% & 1\\
& Taylor mx & 2.1 & 0.672\% & 0\% & 1\\
& Taylor mprev &   2 & 0.672\% & 0\% & 1\\
& Taylor xprev &   2 & 0.672\% & 0\% & 1\\
\hline
\multirow{5}*{$-1$}
& power & 0.16 & 100\% & 100\% & $2^{24}$\\
& explog & 0.18 & 1.27\% & 1.15e-07\% & 2\\
& Taylor mx &   2 & 0.621\% & 0\% & 1\\
& Taylor mprev &   2 & 0.621\% & 0\% & 1\\
& Taylor xprev &   2 & 0.621\% & 0\% & 1\\
\hline
\multirow{5}*{$1/7$}
& power & 0.16 & 100\% & 100\% & ${} > 2^{23}$\\
& explog & 0.18 & 0.81\% & 0\% & 1\\
& Taylor mx &   2 & 0.901\% & 0\% & 1\\
& Taylor mprev &   2 & 0.901\% & 0\% & 1\\
& Taylor xprev &   2 & 0.901\% & 0\% & 1\\
\hline
\multirow{5}*{$1/5$}
& power & 0.16 & 100\% & 100\% & ${} > 2^{23}$\\
& explog & 0.18 & 0.71\% & 0\% & 1\\
& Taylor mx &   2 & 0.854\% & 0\% & 1\\
& Taylor mprev &   2 & 0.854\% & 0\% & 1\\
& Taylor xprev &   2 & 0.854\% & 0\% & 1\\
\hline
\multirow{5}*{$1/3$}
& power & 0.16 & 100\% & 100\% & ${} > 2^{23}$\\
& explog & 0.18 & 0.536\% & 0\% & 1\\
& Taylor mx &   2 & 0.641\% & 0\% & 1\\
& Taylor mprev &   2 & 0.641\% & 0\% & 1\\
& Taylor xprev &   2 & 0.641\% & 0\% & 1\\
\hline
\multirow{5}*{$1$}
& power & 0.16 & 100\% & 100\% & ${} > 2^{23}$\\
& explog & 0.18 & 1.15e-07\% & 0\% & 1\\
& Taylor mx &   2 & 0\% & 0\% & 0\\
& Taylor mprev &   2 & 0\% & 0\% & 0\\
& Taylor xprev &   2 & 0\% & 0\% & 0\\
\hline
\multirow{5}*{$2$}
& power & 0.16 & 100\% & 100\% & ${} > 2^{23}$\\
& explog & 0.18 & 1.27\% & 0\% & 1\\
& Taylor mx &   2 & 0.563\% & 0\% & 1\\
& Taylor mprev &   2 & 0.563\% & 0\% & 1\\
& Taylor xprev &   2 & 0.563\% & 0\% & 1\\
\hline
\multirow{5}*{$3$}
& power & 0.16 & 100\% & 100\% & ${} > 2^{23}$\\
& explog & 0.18 & 7.34\% & 0\% & 1\\
& Taylor mx &   2 & 7.39\% & 0.0372\% & 2\\
& Taylor mprev &   2 & 7.39\% & 0.0372\% & 2\\
& Taylor xprev &   2 & 7.39\% & 0.0372\% & 2\\
\hline
\multirow{5}*{$5$}
& power & 0.16 & 100\% & 100\% & ${} > 2^{23}$\\
& explog & 0.18 & 5.09\% & 0.00161\% & 2\\
& Taylor mx &   2 & 5.25\% & 0.0137\% & 2\\
& Taylor mprev &   2 & 5.25\% & 0.0137\% & 2\\
& Taylor xprev &   2 & 5.25\% & 0.0137\% & 2\\
\hline
\multirow{5}*{$7$}
& power & 0.16 & 100\% & 100\% & ${} > 2^{23}$\\
& explog & 0.18 & 4.16\% & 0\% & 1\\
& Taylor mx &   2 & 4.52\% & 0.0182\% & 2\\
& Taylor mprev &   2 & 4.52\% & 0.0182\% & 2\\
& Taylor xprev &   2 & 4.52\% & 0.0182\% & 2\\
\end{tabular}
\end{table}

\begin{table}
\caption{Performance, in $10^9$ elements (GE) per seconds, and accuracy of the $\powipmi$ function in the power, explog and Taylor forms
in the $[2^{-23}, 2^{-6}[$ range, $1,015,021,568$ single-precision floating-point values.}
\label{tab:pow1pm1.log.3}
\scriptsize\centering
\begin{tabular}{l|l lrrl}
$\alpha$ & method & GE/sec & discrepancies & ${}> 1$ ULP & max ULPs\\
\hline
\multirow{5}*{$-2$}
& power & 0.15 & 100\% & 100\% & $2^{24}$\\
& explog & 0.14 & 3.46\% & 0\% & 1\\
& Taylor mx &   2 & 3.28\% & 0\% & 1\\
& Taylor mprev &   2 & 3.28\% & 0\% & 1\\
& Taylor xprev &   2 & 3.28\% & 0\% & 1\\
\hline
\multirow{5}*{$-1$}
& power & 0.16 & 100\% & 100\% & $2^{24}$\\
& explog & 0.15 & 6.71\% & 1.97e-07\% & 2\\
& Taylor mx &   2 & 3.27\% & 0\% & 1\\
& Taylor mprev &   2 & 3.27\% & 0\% & 1\\
& Taylor xprev &   2 & 3.26\% & 0\% & 1\\
\hline
\multirow{5}*{$1/7$}
& power & 0.16 & 100\% & 100\% & $2^{24}$\\
& explog & 0.15 & 5.46\% & 0.000724\% & 2\\
& Taylor mx &   2 & 4.94\% & 0.000473\% & 2\\
& Taylor mprev &   2 & 4.94\% & 0.000476\% & 2\\
& Taylor xprev & 1.9 & 4.94\% & 0.000475\% & 2\\
\hline
\multirow{5}*{$1/5$}
& power & 0.16 & 100\% & 100\% & ${} > 2^{23}$\\
& explog & 0.15 & 5.41\% & 0.0212\% & 2\\
& Taylor mx &   2 & 4.92\% & 0.000223\% & 2\\
& Taylor mprev &   2 & 4.92\% & 0.000227\% & 2\\
& Taylor xprev & 1.9 & 4.92\% & 0.000225\% & 2\\
\hline
\multirow{5}*{$1/3$}
& power & 0.16 & 100\% & 100\% & ${} > 2^{23}$\\
& explog & 0.15 & 5.35\% & 0.00742\% & 2\\
& Taylor mx &   2 & 4.85\% & 8.73e-05\% & 2\\
& Taylor mprev &   2 & 4.85\% & 8.79e-05\% & 2\\
& Taylor xprev &   2 & 4.85\% & 8.73e-05\% & 2\\
\hline
\multirow{5}*{$1$}
& power & 0.16 & 100\% & 100\% & ${} > 2^{23}$\\
& explog & 0.15 & 0.0203\% & 0\% & 1\\
& Taylor mx &   2 & 0\% & 0\% & 0\\
& Taylor mprev &   2 & 0\% & 0\% & 0\\
& Taylor xprev &   2 & 0\% & 0\% & 0\\
\hline
\multirow{5}*{$2$}
& power & 0.16 & 100\% & 100\% & ${} > 2^{23}$\\
& explog & 0.15 & 6.72\% & 0\% & 1\\
& Taylor mx &   2 & 5.8\% & 0\% & 1\\
& Taylor mprev &   2 & 5.8\% & 0\% & 1\\
& Taylor xprev &   2 & 5.8\% & 0\% & 1\\
\hline
\multirow{5}*{$3$}
& power & 0.16 & 100\% & 100\% & ${} > 2^{23}$\\
& explog & 0.14 & 11.4\% & 0.067\% & 2\\
& Taylor mx &   2 & 12.2\% & 0.142\% & 2\\
& Taylor mprev &   2 & 12.2\% & 0.142\% & 2\\
& Taylor xprev &   2 & 12.2\% & 0.142\% & 2\\
\hline
\multirow{5}*{$5$}
& power & 0.16 & 100\% & 100\% & ${} > 2^{23}$\\
& explog & 0.14 & 9.25\% & 0.0176\% & 2\\
& Taylor mx &   2 & 10.3\% & 0.0925\% & 2\\
& Taylor mprev &   2 & 10.3\% & 0.0925\% & 2\\
& Taylor xprev &   2 & 10.3\% & 0.0926\% & 2\\
\hline
\multirow{5}*{$7$}
& power & 0.16 & 100\% & 100\% & ${} > 2^{23}$\\
& explog & 0.14 & 8.4\% & 0.029\% & 2\\
& Taylor mx &   2 & 9.73\% & 0.0905\% & 2\\
& Taylor mprev &   2 & 9.73\% & 0.0905\% & 2\\
& Taylor xprev &   2 & 9.73\% & 0.0905\% & 2\\
\end{tabular}
\end{table}

\begin{table}
\caption{Performance, in $10^9$ elements (GE) per seconds, and accuracy of the $\powipmi$ function in the power, explog and Taylor forms
in the $[2^{-6}, 2^{-3}[$ range, covering $1,040,187,392$ single-precision floating-point values.}
\label{tab:pow1pm1.log.4}
\scriptsize\centering
\begin{tabular}{l|l lrrl}
$\alpha$ & method & GE/sec & discrepancies & ${}> 1$ ULP & max ULPs\\
\hline
\multirow{5}*{$-2$}
& power & 0.15 & 99.9\% & 99.6\% & $2^{24}$\\
& explog & 0.14 & 3.95\% & 6.73e-07\% & 2\\
& Taylor mx & 1.8 & 3.94\% & 0.187\% & ${} > 2^{2}$\\
& Taylor mprev & 1.8 & 3.94\% & 0.187\% & ${} > 2^{2}$\\
& Taylor xprev & 1.8 & 3.94\% & 0.187\% & ${} > 2^{2}$\\
\hline
\multirow{5}*{$-1$}
& power & 0.16 & 99.9\% & 99.7\% & $2^{24}$\\
& explog & 0.14 & 7.13\% & 6.73e-07\% & 2\\
& Taylor mx & 1.8 & 3.72\% & 0.00255\% & 2\\
& Taylor mprev & 1.8 & 3.72\% & 0.00255\% & 2\\
& Taylor xprev & 1.8 & 3.71\% & 0.00175\% & 2\\
\hline
\multirow{5}*{$1/7$}
& power & 0.16 & 100\% & 100\% & $2^{24}$\\
& explog & 0.14 & 6.15\% & 0.00104\% & 2\\
& Taylor mx & 1.8 & 5.55\% & 0.00245\% & 2\\
& Taylor mprev & 1.7 & 5.55\% & 0.00256\% & 2\\
& Taylor xprev & 1.7 & 5.55\% & 0.00256\% & 2\\
\hline
\multirow{5}*{$1/5$}
& power & 0.16 & 100\% & 99.9\% & ${} > 2^{23}$\\
& explog & 0.14 & 6.1\% & 0.0242\% & 2\\
& Taylor mx & 1.8 & 5.5\% & 0.00122\% & 2\\
& Taylor mprev & 1.7 & 5.51\% & 0.00134\% & 2\\
& Taylor xprev & 1.7 & 5.51\% & 0.0013\% & 2\\
\hline
\multirow{5}*{$1/3$}
& power & 0.16 & 100\% & 99.9\% & ${} > 2^{23}$\\
& explog & 0.14 & 6.05\% & 0.00914\% & 2\\
& Taylor mx & 1.8 & 5.44\% & 0.00121\% & 2\\
& Taylor mprev & 1.7 & 5.45\% & 0.00128\% & 2\\
& Taylor xprev & 1.7 & 5.45\% & 0.00122\% & 2\\
\hline
\multirow{5}*{$1$}
& power & 0.16 & 99.9\% & 99.7\% & ${} > 2^{23}$\\
& explog & 0.14 & 0.149\% & 0\% & 1\\
& Taylor mx & 1.8 & 0\% & 0\% & 0\\
& Taylor mprev & 1.8 & 0\% & 0\% & 0\\
& Taylor xprev & 1.8 & 0\% & 0\% & 0\\
\hline
\multirow{5}*{$2$}
& power & 0.16 & 99.9\% & 99.7\% & ${} > 2^{23}$\\
& explog & 0.14 & 7.2\% & 0\% & 1\\
& Taylor mx & 1.8 & 6.52\% & 0\% & 1\\
& Taylor mprev & 1.8 & 6.52\% & 0\% & 1\\
& Taylor xprev & 1.8 & 6.52\% & 0\% & 1\\
\hline
\multirow{5}*{$3$}
& power & 0.16 & 99.9\% & 99.7\% & ${} > 2^{23}$\\
& explog & 0.14 & 12\% & 0.0834\% & 2\\
& Taylor mx & 1.8 & 12.9\% & 0.157\% & 2\\
& Taylor mprev & 1.8 & 12.9\% & 0.157\% & 2\\
& Taylor xprev & 1.8 & 12.9\% & 0.157\% & 2\\
\hline
\multirow{5}*{$5$}
& power & 0.16 & 99.9\% & 99.7\% & ${} > 2^{23}$\\
& explog & 0.14 & 10\% & 0.0412\% & 2\\
& Taylor mx & 1.8 & 11\% & 0.102\% & 2\\
& Taylor mprev & 1.8 & 11\% & 0.102\% & 2\\
& Taylor xprev & 1.7 & 11\% & 0.103\% & 2\\
\hline
\multirow{5}*{$7$}
& power & 0.16 & 99.9\% & 99.7\% & ${} > 2^{23}$\\
& explog & 0.14 & 9.19\% & 0.0576\% & 2\\
& Taylor mx & 1.8 & 10.4\% & 0.1\% & 2\\
& Taylor mprev & 1.8 & 10.4\% & 0.1\% & 2\\
& Taylor xprev & 1.8 & 10.4\% & 0.1\% & 2\\
\end{tabular}
\end{table}

\begin{table}
\caption{Performance, in $10^9$ elements (GE) per seconds, and accuracy of the $\powipmi$ function in the power, explog and Taylor forms
in the $[2^{-3}, 1[$ range, covering $1,065,353,216$ single-precision floating-point values.}
\label{tab:pow1pm1.log.5}
\scriptsize\centering
\begin{tabular}{l|l lrrl}
$\alpha$ & method & GE/sec & discrepancies & ${}> 1$ ULP & max ULPs\\
\hline
\multirow{5}*{$-2$}
& power & 0.15 & 98.9\% & 97.8\% & $2^{24}$\\
& explog & 0.14 & 4.39\% & 9.1e-06\% & 2\\
& Taylor mx &   2 & 6.21\% & 2.54\% & ${} > 2^{26}$\\
& Taylor mprev &   2 & 6.21\% & 2.54\% & ${} > 2^{26}$\\
& Taylor xprev &   2 & 6.21\% & 2.54\% & ${} > 2^{26}$\\
\hline
\multirow{5}*{$-1$}
& power & 0.16 & 99.1\% & 97.9\% & $2^{24}$\\
& explog & 0.14 & 7.57\% & 9.88e-05\% & 2\\
& Taylor mx &   2 & 6\% & 2.31\% & ${} > 2^{23}$\\
& Taylor mprev &   2 & 6\% & 2.31\% & ${} > 2^{23}$\\
& Taylor xprev &   2 & 5.98\% & 2.31\% & ${} > 2^{23}$\\
\hline
\multirow{5}*{$1/7$}
& power & 0.16 & 99.9\% & 99.7\% & $2^{24}$\\
& explog & 0.14 & 6.87\% & 0.00435\% & 2\\
& Taylor mx &   2 & 7.58\% & 1.97\% & ${} > 2^{19}$\\
& Taylor mprev &   2 & 7.58\% & 1.97\% & ${} > 2^{19}$\\
& Taylor xprev &   2 & 7.58\% & 1.97\% & ${} > 2^{19}$\\
\hline
\multirow{5}*{$1/5$}
& power & 0.16 & 99.9\% & 99.6\% & ${} > 2^{23}$\\
& explog & 0.14 & 6.83\% & 0.0344\% & 2\\
& Taylor mx &   2 & 7.49\% & 1.89\% & ${} > 2^{18}$\\
& Taylor mprev &   2 & 7.51\% & 1.89\% & ${} > 2^{18}$\\
& Taylor xprev &   2 & 7.51\% & 1.9\% & ${} > 2^{18}$\\
\hline
\multirow{5}*{$1/3$}
& power & 0.16 & 99.8\% & 99.3\% & ${} > 2^{23}$\\
& explog & 0.14 & 6.82\% & 0.0223\% & 2\\
& Taylor mx &   2 & 7.41\% & 1.82\% & ${} > 2^{18}$\\
& Taylor mprev &   2 & 7.42\% & 1.82\% & ${} > 2^{18}$\\
& Taylor xprev &   2 & 7.42\% & 1.82\% & ${} > 2^{18}$\\
\hline
\multirow{5}*{$1$}
& power & 0.16 & 99.2\% & 98\% & ${} > 2^{23}$\\
& explog & 0.14 & 0.758\% & 9.39e-08\% & 2\\
& Taylor mx &   2 & 0\% & 0\% & 0\\
& Taylor mprev &   2 & 0\% & 0\% & 0\\
& Taylor xprev &   2 & 0\% & 0\% & 0\\
\hline
\multirow{5}*{$2$}
& power & 0.16 & 99.3\% & 98.2\% & ${} > 2^{23}$\\
& explog & 0.13 & 7.94\% & 0.0205\% & ${} > 2$\\
& Taylor mx &   2 & 7.07\% & 0\% & 1\\
& Taylor mprev &   2 & 7.07\% & 0\% & 1\\
& Taylor xprev &   2 & 7.07\% & 0\% & 1\\
\hline
\multirow{5}*{$3$}
& power & 0.16 & 99.3\% & 98.4\% & ${} > 2^{23}$\\
& explog & 0.13 & 13.1\% & 0.273\% & $2^{2}$\\
& Taylor mx &   2 & 13.5\% & 0.174\% & 2\\
& Taylor mprev &   2 & 13.5\% & 0.17\% & 2\\
& Taylor xprev &   2 & 13.5\% & 0.182\% & ${} > 2$\\
\hline
\multirow{5}*{$5$}
& power & 0.16 & 99.3\% & 98.7\% & ${} > 2^{23}$\\
& explog & 0.13 & 11.2\% & 0.423\% & ${} > 2^{2}$\\
& Taylor mx &   2 & 11.7\% & 0.122\% & ${} > 2$\\
& Taylor mprev &   2 & 11.7\% & 0.122\% & ${} > 2$\\
& Taylor xprev &   2 & 11.8\% & 0.151\% & $2^{2}$\\
\hline
\multirow{5}*{$7$}
& power & 0.16 & 99.3\% & 98.9\% & ${} > 2^{23}$\\
& explog & 0.13 & 10.6\% & 0.639\% & ${} > 2^{3}$\\
& Taylor mx &   2 & 11.3\% & 0.167\% & $2^{2}$\\
& Taylor mprev &   2 & 11.2\% & 0.166\% & $2^{2}$\\
& Taylor xprev &   2 & 11.2\% & 0.174\% & ${} > 2^{2}$\\
\end{tabular}
\end{table}

Results for the CPU benchmark are collected in Tables~\ref{tab:pow1pm1.log.1}--\ref{tab:pow1pm1.log.5}.
Of note is the marginal performance advantage of the power form over the explog form in the higher ranges,
at the expense of an (expected) massive drop in accuracy in the entire domain,
and the significant performance advantage of the Taylor forms, that approaches the
memory bandwidth limit with a numerical accuracy that is comparable with,
and sometimes better than, that of the explog form,
with a slight advantage for the `mx' variant both in performance and accuracy.
The explog form turns out to be exact in the denormals, and maintains a maximum ULP of 2 until the ``perten''
range (argument larger than $1/8$), and even then it only gets worse results for less than $1\%$ of the values.

Aside from the computational performance advantage,
that in this case is skewed by the cross-iteration vectorization,
a significant advantage of the Taylor forms is the correctly rounded results given for $\alpha = 1$:
adopting this form gives the best results for the common case $\gamma = 1$ even without special-casing
the choice of polytropic constant;
even though lower orders would also give the same results with some further performance gains,
adopting the Taylor expansion at order 8 allows keeping the code generic without particular losses
in either performance or accuracy regardless of the choice of $\gamma$.

In fact, the only range\slash exponent combination in which the Taylor expansion falls significantly short
of the explog form in terms of accuracy is in the ``perten'' range for $\alpha < 1$ (Table~\ref{tab:pow1pm1.log.5}).
In this case, the accuracy achieved with an expansion of order 8 quickly drops off as $x \to 1$.
Such an application of the $\powipmi$ function is not expected in the normal applications of WCSPH,
but it may be safer to adopt the explog for the inverse of the equation of state
for larger values of $\gamma$.

\subsection{Equation of state round-trip}\label{sec:roundtrip}

A second but important measure
in the evaluation of the best choice for the expression of the $\powipmi$ function
is the accuracy of the `round-trip',
i.e. the accuracy with which one can recover the original (scaled) density from the computed (scaled)
pressure by applying the inverse equation of state, or conversely:
considering that analytically we should have $\powipmi( \powipmi(x, 1/\alpha), \alpha) = x$,
we can evaluate the round-trip accuracy of the $\powipmi$ forms for different exponents
by comparing the result of the chain evaluation of the function and its inverse against
the original argument. Once again, accuracy can be estimated in terms of number of elements
for which the results are different, number of elements for which the maximum error is
larger than 1 ULP, and the maximum ULPs in the error.

\begin{table}
\caption{Round-trip accuracy of the $\powipmi$ function in the explog and Taylor (mx) forms.
This table refers to the range $[0, 2^{-126}[$ that contains $8,388,608$ single-precision floating-point values.}
\label{tab:pow1pm1.rt.log.1}
\scriptsize\centering
\begin{tabular}{l|l rrl}
$\alpha$ & method & discrepancies & ${}> 1$ ULP & max ULPs\\
\hline
\multirow{2}*{$-2$}
& explog & 0\% & 0\% & 0\\
& Taylor mx & 0\% & 0\% & 0\\
\hline
\multirow{2}*{$-1$}
& explog & 0\% & 0\% & 0\\
& Taylor mx & 0\% & 0\% & 0\\
\hline
\multirow{2}*{$1/7$}
& explog & 85.7\% & 57.1\% & $2^{2}$\\
& Taylor mx & 85.7\% & 57.1\% & $2^{2}$\\
\hline
\multirow{2}*{$1/5$}
& explog & 80\% & 40\% & 2\\
& Taylor mx & 80\% & 40\% & 2\\
\hline
\multirow{2}*{$1/3$}
& explog & 66.7\% & 0\% & 1\\
& Taylor mx & 66.7\% & 0\% & 1\\
\hline
\multirow{2}*{$1$}
& explog & 0\% & 0\% & 0\\
& Taylor mx & 0\% & 0\% & 0\\
\hline
\multirow{2}*{$2$}
& explog & 0\% & 0\% & 0\\
& Taylor mx & 0\% & 0\% & 0\\
\hline
\multirow{2}*{$3$}
& explog & 8.33\% & 0\% & 1\\
& Taylor mx & 8.33\% & 0\% & 1\\
\hline
\multirow{2}*{$5$}
& explog & 2.5\% & 0\% & 1\\
& Taylor mx & 2.5\% & 0\% & 1\\
\hline
\multirow{2}*{$7$}
& explog & 6.55\% & 0\% & 1\\
& Taylor mx & 6.55\% & 0\% & 1\\
\end{tabular}
\end{table}

\begin{table}
\caption{Round-trip accuracy of the $\powipmi$ function in the explog and Taylor (mx) forms.
This table refers to the range $[2^{-126}, 2^{-23}[$ that contains $872,415,232$ single-precision floating-point values.}
\label{tab:pow1pm1.rt.log.2}
\scriptsize\centering
\begin{tabular}{l|l rrl}
$\alpha$ & method & discrepancies & ${}> 1$ ULP & max ULPs\\
\hline
\multirow{2}*{$-2$}
& explog & 3.44e-07\% & 0\% & 1\\
& Taylor mx & 1.76\% & 0\% & 1\\
\hline
\multirow{2}*{$-1$}
& explog & 2.29e-07\% & 1.15e-07\% & 2\\
& Taylor mx & 0.962\% & 0\% & 1\\
\hline
\multirow{2}*{$1/7$}
& explog & 29.9\% & 1.27\% & $2^{2}$\\
& Taylor mx & 30.5\% & 1.53\% & $2^{2}$\\
\hline
\multirow{2}*{$1/5$}
& explog & 24\% & 0.769\% & 2\\
& Taylor mx & 24.5\% & 0.841\% & 2\\
\hline
\multirow{2}*{$1/3$}
& explog & 34.1\% & 0\% & 1\\
& Taylor mx & 34.3\% & 0.0801\% & 2\\
\hline
\multirow{2}*{$1$}
& explog & 1.15e-07\% & 0\% & 1\\
& Taylor mx & 0\% & 0\% & 0\\
\hline
\multirow{2}*{$2$}
& explog & 2.29e-07\% & 0\% & 1\\
& Taylor mx & 0.962\% & 0\% & 1\\
\hline
\multirow{2}*{$3$}
& explog & 33.1\% & 0\% & 1\\
& Taylor mx & 33.6\% & 0\% & 1\\
\hline
\multirow{2}*{$5$}
& explog & 22.3\% & 0\% & 1\\
& Taylor mx & 23.2\% & 0.0881\% & 2\\
\hline
\multirow{2}*{$7$}
& explog & 54.3\% & 2.29e-07\% & 2\\
& Taylor mx & 54.3\% & 0.381\% & ${} > 2$\\
\end{tabular}
\end{table}

\begin{table}
\caption{Round-trip accuracy of the $\powipmi$ function in the explog and Taylor (mx) forms.
This table refers to the range $[2^{-23}, 2^{-6}[$ that contains $1,015,021,568$ single-precision floating-point values.}
\label{tab:pow1pm1.rt.log.3}
\scriptsize\centering
\begin{tabular}{l|l rrl}
$\alpha$ & method & discrepancies & ${}> 1$ ULP & max ULPs\\
\hline
\multirow{2}*{$-2$}
& explog & 0.0602\% & 0.000435\% & 2\\
& Taylor mx & 7.02\% & 0.0291\% & 2\\
\hline
\multirow{2}*{$-1$}
& explog & 0.0403\% & 0.000229\% & 2\\
& Taylor mx & 6.38\% & 0.0282\% & 2\\
\hline
\multirow{2}*{$1/7$}
& explog & 29.7\% & 1.1\% & $2^{2}$\\
& Taylor mx & 33.8\% & 2.18\% & $2^{2}$\\
\hline
\multirow{2}*{$1/5$}
& explog & 23.8\% & 0.668\% & ${} > 2$\\
& Taylor mx & 28.3\% & 1.54\% & ${} > 2$\\
\hline
\multirow{2}*{$1/3$}
& explog & 34\% & 0.00955\% & ${} > 2$\\
& Taylor mx & 37.1\% & 1.35\% & $2^{2}$\\
\hline
\multirow{2}*{$1$}
& explog & 0.0203\% & 0.00108\% & 2\\
& Taylor mx & 0\% & 0\% & 0\\
\hline
\multirow{2}*{$2$}
& explog & 0.0364\% & 0.000214\% & 2\\
& Taylor mx & 6.05\% & 0.0254\% & 2\\
\hline
\multirow{2}*{$3$}
& explog & 33.1\% & 0.00961\% & ${} > 2$\\
& Taylor mx & 36.5\% & 0.981\% & $2^{2}$\\
\hline
\multirow{2}*{$5$}
& explog & 22.4\% & 0.00867\% & ${} > 2$\\
& Taylor mx & 27.3\% & 0.927\% & $2^{2}$\\
\hline
\multirow{2}*{$7$}
& explog & 54.4\% & 0.019\% & ${} > 2$\\
& Taylor mx & 55.1\% & 1.9\% & $2^{2}$\\
\end{tabular}
\end{table}

\begin{table}
\caption{Round-trip accuracy of the $\powipmi$ function in the explog and Taylor (mx) forms.
This table refers to the range $[2^{-6}, 2^{-3}[$ that contains $1,040,187,392$ single-precision floating-point values.}
\label{tab:pow1pm1.rt.log.4}
\scriptsize\centering
\begin{tabular}{l|l rrl}
$\alpha$ & method & discrepancies & ${}> 1$ ULP & max ULPs\\
\hline
\multirow{2}*{$-2$}
& explog & 0.413\% & 0.0225\% & 2\\
& Taylor mx & 8.16\% & 0.547\% & ${} > 2^{5}$\\
\hline
\multirow{2}*{$-1$}
& explog & 0.287\% & 0.0118\% & 2\\
& Taylor mx & 7.25\% & 0.0845\% & ${} > 2$\\
\hline
\multirow{2}*{$1/7$}
& explog & 29.7\% & 1.13\% & $2^{2}$\\
& Taylor mx & 34.3\% & 2.28\% & $2^{2}$\\
\hline
\multirow{2}*{$1/5$}
& explog & 23.9\% & 0.693\% & ${} > 2$\\
& Taylor mx & 28.8\% & 1.63\% & ${} > 2$\\
\hline
\multirow{2}*{$1/3$}
& explog & 34.1\% & 0.0675\% & ${} > 2$\\
& Taylor mx & 37.6\% & 1.53\% & $2^{2}$\\
\hline
\multirow{2}*{$1$}
& explog & 0.149\% & 0.00815\% & 2\\
& Taylor mx & 0\% & 0\% & 0\\
\hline
\multirow{2}*{$2$}
& explog & 0.253\% & 0.00913\% & 2\\
& Taylor mx & 7.04\% & 0.305\% & ${} > 2^{3}$\\
\hline
\multirow{2}*{$3$}
& explog & 33.3\% & 0.0742\% & ${} > 2$\\
& Taylor mx & 37.3\% & 1.79\% & ${} > 2^{9}$\\
\hline
\multirow{2}*{$5$}
& explog & 22.6\% & 0.0732\% & ${} > 2$\\
& Taylor mx & 28.6\% & 2.39\% & ${} > 2^{17}$\\
\hline
\multirow{2}*{$7$}
& explog & 54.4\% & 0.139\% & ${} > 2$\\
& Taylor mx & 56\% & 3.75\% & ${} > 2^{22}$\\
\end{tabular}
\end{table}

\begin{table}
\caption{Round-trip accuracy of the $\powipmi$ function in the explog and Taylor (mx) forms.
This table refers to the range $[2^{-3}, 1[$ that contains $1,065,353,216$ single-precision floating-point values.}
\label{tab:pow1pm1.rt.log.5}
\scriptsize\centering
\begin{tabular}{l|l rrl}
$\alpha$ & method & discrepancies & ${}> 1$ ULP & max ULPs\\
\hline
\multirow{2}*{$-2$}
& explog & 1.55\% & 0.307\% & ${} > 2$\\
& Taylor mx & 10.3\% & 2.9\% & ${} > 2^{37}$\\
\hline
\multirow{2}*{$-1$}
& explog & 1.2\% & 0.13\% & ${} > 2$\\
& Taylor mx & 9.43\% & 2.38\% & ${} > 2^{23}$\\
\hline
\multirow{2}*{$1/7$}
& explog & 30.1\% & 1.29\% & ${} > 2^{2}$\\
& Taylor mx & 35.7\% & 4.24\% & ${} > 2^{20}$\\
\hline
\multirow{2}*{$1/5$}
& explog & 24.3\% & 0.843\% & ${} > 2^{2}$\\
& Taylor mx & 30.3\% & 3.53\% & ${} > 2^{19}$\\
\hline
\multirow{2}*{$1/3$}
& explog & 34.4\% & 0.263\% & $2^{2}$\\
& Taylor mx & 38.8\% & 3.33\% & ${} > 2^{19}$\\
\hline
\multirow{2}*{$1$}
& explog & 0.758\% & 0.0315\% & 2\\
& Taylor mx & 0\% & 0\% & 0\\
\hline
\multirow{2}*{$2$}
& explog & 0.965\% & 0.0285\% & ${} > 2$\\
& Taylor mx & 9.23\% & 2.66\% & ${} > 2^{29}$\\
\hline
\multirow{2}*{$3$}
& explog & 33.6\% & 0.312\% & $2^{2}$\\
& Taylor mx & 38.8\% & 4.11\% & ${} > 2^{40}$\\
\hline
\multirow{2}*{$5$}
& explog & 23\% & 0.222\% & ${} > 2$\\
& Taylor mx & 30.3\% & 4.7\% & ${} > 2^{57}$\\
\hline
\multirow{2}*{$7$}
& explog & 54.6\% & 0.532\% & ${} > 2$\\
& Taylor mx & 57\% & 6.03\% & ${} > 2^{73}$\\
\end{tabular}
\end{table}

We restrict our analysis to the explog and ``Taylor mx'' forms,
taken (from our previous analysis) respectively
as the most accurate in the full range,
and the fastest.
The results are presented in Tables~\ref{tab:pow1pm1.rt.log.1} to~\ref{tab:pow1pm1.rt.log.5}.

It's interesting to note that the round-trip accuracy is not the same when starting from $\alpha$ or from $1/\alpha$:
the error distribution for $\powipmi( \powipmi(x, \alpha), 1/\alpha)$ is different from that of $\powipmi( \powipmi(x, 1/\alpha), \alpha)$.
This can be observed for the exponents $3, 5, 7$ for which we also have the results of the round-trip
starting from the reciprocals $1/3, 1/5, 1/7$.

In the denormals range in particular (Table~\ref{tab:pow1pm1.rt.log.1}) this affects a majority of the values,
due to the fact that starting from the fractional ($0 < \alpha < 1$) exponent leads to loss of bits in the mantissa,
a loss that is not experienced when starting from exponents $\alpha \ge 1$.

Despite this asymmetry, the round-trip error is limited to 2~significant bits at most for either form until the percent range,
where the Taylor form starts losing accuracy for exponents different from~$1$,
with the effect becoming catastrophic in the perten range:
this is expected, as we have already seen that in this range the Taylor form has low accuracy for $\alpha < 1$,
and in the round-trip tests either $\alpha$ or its reciprocal will be less than~1
whenever $\alpha \ne 1$.
On the upside, the Taylor form gives exact results throughout the whole range for $\alpha = 1$,
whereas explog loses one or two ULPs in the upper subranges.

\section{Applications to the ``dummy'' boundary model}\label{sec:dummy}

Accurate evaluation of the pressure is essential in all formulations of WCSPH,
but it becomes even more important in boundary models where its evaluation is needed for the correct imposition of boundary conditions.

As discussed in Section~\ref{sec:dummy_intro}, this is the case for the dummy boundary model.
In this case, it's important not only that the pressure evaluation is accurate,
but also that the round-trip through the inverse of the equation of state does not introduce significant perturbations.
Although the formula we present here are specifically intended for the dummy boundary model,
similar formulas can be derived for any solid boundary model
(or in fact any other context)
that requires the computation of the pressure from neighboring particles.

Let us consider at first the pressure interpolation (equation~\eqref{eq:dummy_pressure})
and rewrite it for the scaled pressure, assuming at first that we only have a single fluid
with bulk modulus $B = \rho_0 c_0^2/\gamma$. Considering that $\rho_j = \rho_0(1+\tilde\rho_j)$, we have
\[
\tilde P_w = \frac{
\sum_j (\tilde P_j + (\gamma/c_0^2) (1+\tilde\rho_j) \vec g \cdot \vec x_{wj}) W(\vec x_{wj}, h)
}{
\sum_j W(\vec x_{wj}, h)
}.
\]
Given the potentially different magnitude of the terms, it is better to explicitly split
the expression in the interpolating and hydrostatic contributions:
\begin{equation}
\tilde P_w = \frac{
\sum_j \tilde P_j W(\vec x_{wj}, h)
}{
\sum_j W(\vec x_{wj}, h)
}
+
\frac{\gamma}{c_0^2}
\frac{
\sum_j (1+\tilde\rho_j) \vec g \cdot \vec x_{wj} W(\vec x_{wj}, h)
}{
\sum_j W(\vec x_{wj}, h)
}.\label{eq:dummy_tilde}
\end{equation}

If multiple fluids are present, each with their bulk module, reference density, speed of sound,
some simplifications will not be possible, and should be replaced with appropriate ratios;
we have $P_j/B_w = P_j/B_j \cdot B_j/B_w = \tilde P_j \cdot \tilde B_{j/w}$ with
$\tilde B_{j/w} = B_j/B_w$ for the pressure interpolation terms, and
$\rho_j/B_w = (\gamma/c_{0,w}^2)\rho_{0,j}(1 + \tilde\rho_j)/\rho_{0,w} =
(\gamma/c_{0,w}^2)\tilde\rho_{0,j/w}(1 + \tilde\rho_j)$ with
$\tilde\rho_{0,j/w} = \rho_{0,j}/\rho_{0,w}$ for the hydrostatic term,
so that equation~\eqref{eq:dummy_tilde} can be rewritten as:
\begin{equation}
\tilde P_w = \frac{
\sum_j \tilde B_{j/w} \tilde P_j W(\vec x_{wj}, h)
}{
\sum_j W(\vec x_{wj}, h)
}
+
\frac{\gamma}{c_{0,w}^2}
\frac{
\sum_j \tilde \rho_{0,j/w} (1+\tilde\rho_j) \vec g \cdot \vec x_{wj} W(\vec x_{wj}, h)
}{
\sum_j W(\vec x_{wj}, h)
}.\label{eq:dummy_tilde_multi}
\end{equation}
in the multi-fluid case.

The wall particle pressure computed as \eqref{eq:dummy_tilde} or \eqref{eq:dummy_tilde_multi}
is used both to compute the wall boundary contribution to the pressure gradient of neighboring fluid particles,
and to compute the density of the wall boundary particles using the inverse equation of state.

In fact, computing the density is important for both uses, since the pressure gradient contribution,
is frequently computed as $P/\rho^k$ for some exponent $k$,
with $k=2$ being a common exponent, such as in the formulation we are using (equation~\eqref{eq:ns-sph1}).

\subsection{Pressure gradient contribution}\label{sec:dummy-pressure}

We have, for any exponent $k$ and for any particle:
\begin{equation}
\frac{P}{\rho^k} = \frac{ B \tilde P}{\rho_0^k (1 + \tilde \rho)^k} =
\frac{c_0^2}{\gamma \rho_0^{k-1}} \frac{ \tilde P }{ (1 + \tilde\rho)^k }.\label{eq:Prhok}
\end{equation}
For fluid particles,
where the scaled pressure is computed directly from the scaled density using the equation of state,
we have:
\begin{align*}
\frac{ \tilde P }{ (1 + \tilde\rho)^k } & {} =
\frac{ (1 + \tilde\rho)^\gamma - 1 }{ (1 + \tilde\rho)^k } =
(1 + \tilde\rho)^{\gamma - k} - (1 + \tilde\rho)^{-k} = {} \\
& {} = \left((1 + \tilde\rho)^{\gamma - k} - 1\right) - \left((1 + \tilde\rho)^{-k} - 1 \right) =
\powipmi(\tilde\rho, \gamma - k) - \powipmi(\tilde\rho, -k)
\end{align*}
so that we can express \eqref{eq:Prhok} entirely in terms of $\powipmi$ times a constant that can be precomputed:
\[
\frac{P}{\rho^k} =
\frac{c_0^2}{\gamma \rho_0^{k-1}} \left(\powipmi(\tilde\rho, \gamma - k) - \powipmi(\tilde\rho, -k)\right)
\]

For dummy boundary particles,
in combination with \eqref{eq:dummy_tilde},
equation~\eqref{eq:Prhok} allows some additional simplifications in the hydrostatic contribution:
\begin{equation}
\frac{P_w}{\rho_w^k} =
\frac{c_{0,w}^2}{\gamma \rho_{0,w}^{k-1}}
\frac{
\sum_j \tilde P_j W(\vec x_{wj}, h)
}{
(1 + \tilde\rho_w)^k \sum_j W(\vec x_{wj}, h)
}
+
\frac{
\sum_j (1+\tilde\rho_j) \vec g \cdot \vec x_{wj} W(\vec x_{wj}, h)
}{
\rho_{0,w}^{k-1}
(1 + \tilde\rho_w)^k \sum_j W(\vec x_{wj}, h)
}.\label{eq:dummy_Prhok_tilde}
\end{equation}
and a similar expression can also be derived for \eqref{eq:dummy_tilde_multi} in the multi-fluid case.
Note that despite the simplifications, \eqref{eq:dummy_Prhok_tilde} still depends on equation~\eqref{eq:dummy_tilde}
to obtain $\tilde\rho_w$ from $\tilde P_w$.
It is however possible to compute the common sub-expressions once and then produce from them
both $\tilde P_w$ for $\tilde\rho_w$ and $\tilde P_w/\rho_w^k$ as needed.

The numerical results we will present momentarily also show that simply using
the expressions~\eqref{eq:dummy_tilde}--\eqref{eq:dummy_Prhok_tilde} over the naive ones
(equation \eqref{eq:dummy_pressure} for the wall boundary, and $P/\rho^k$)
is sufficient to improve numerical accuracy in the computation of the acoustic terms,
regardless of the choice of $\powipmi$ discussed in section~\ref{sec:powforms}.

\subsection{Numerical results}


\paragraph{SPHERIC Test Case 2}

The first test case is the classic SPHERIC~2 reference case,
a classic bore-in-a-box test on a domain which is $3.22\times1\times1 \text{m}^3$,
with a $0.55$m water height behind a closed gate positioned $2$m from the back wall.
A box (dimensions $0.161\times0.403\times0.161$)
is placed with its water-facing side at $1.18$m from the gate,
centered with respect to the side walls.
The box is covered with 8 pressure sensors along the central meridian.
Water gages take the water elevation over time at three locations
(halfway between the gate and the obstacle, near the front of the obstacle,
near the back of the obstacle).
Further details about the setup can be found in~\cite{spheric2_kleefsman,spheric2_lee}.

\begin{figure}[btp]
\includegraphics[width=\columnwidth]{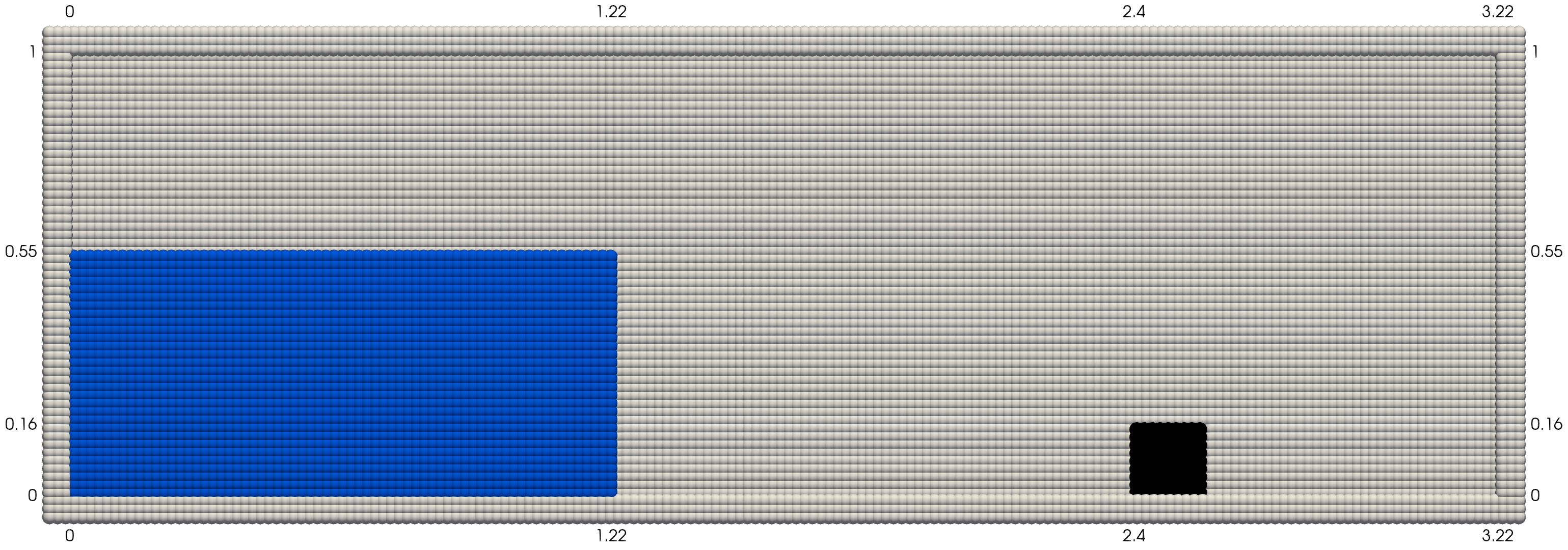}
\caption{Cross-section of the 3D geometry for the SPHERIC~2 Test Case.
Fluid particles are blue, solid boundary particles realizing the walls are gray,
and solid boundary particles realizing the obstacle are black.}
\label{fig:spheric2-t00}
\end{figure}

\begin{figure}
\includegraphics[width=\columnwidth]{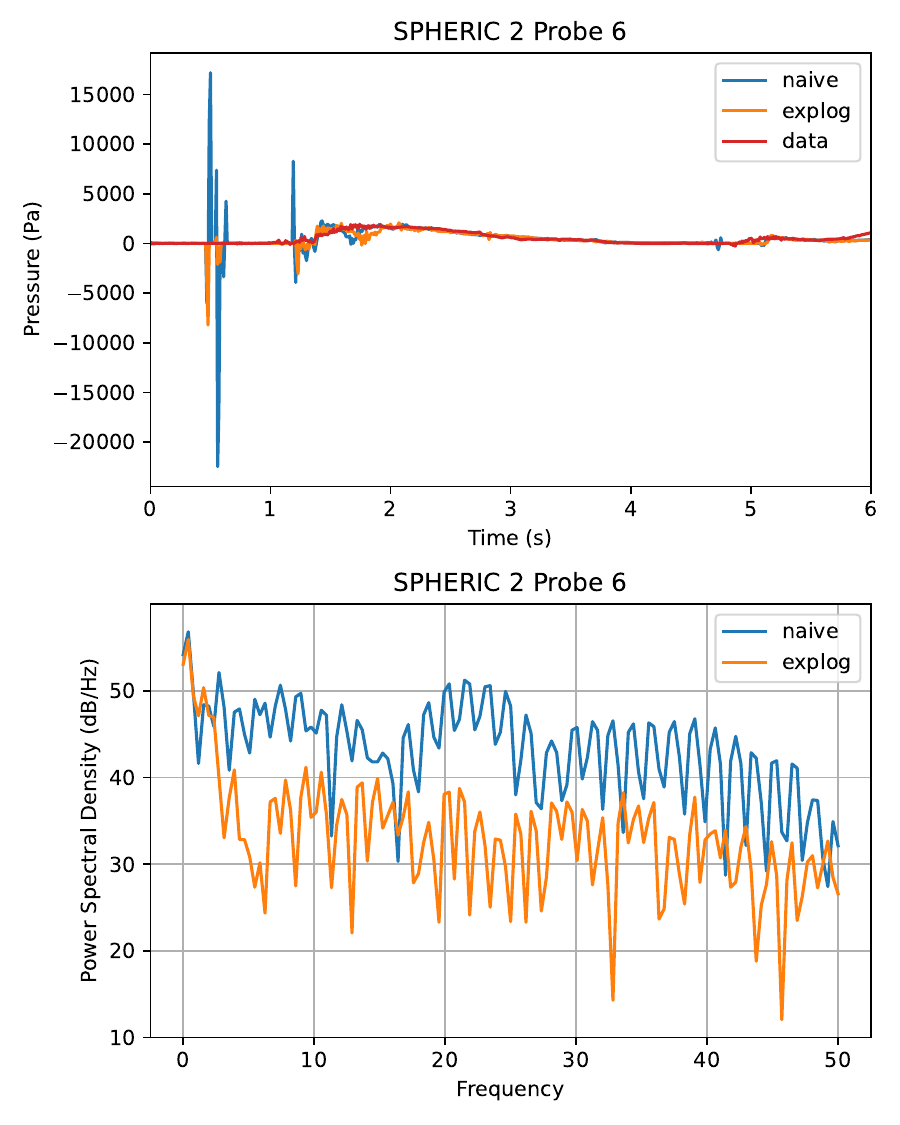}
\caption{Pressure oscillations and their spectral density at Probe~6 of the SPHERIC 2 Test case,
when using the naive versus the explog forms of the acoustic terms.}
\label{fig:fft}
\end{figure}

We model the test case using the dummy solid boundary formulation
(Figure~\ref{fig:spheric2-t00}).
The equation of state for the fluid uses the polytropic constant $\gamma = 7$,
and the Monaghan artificial viscosity is used with $\alpha = 0.05$.
The particle resolution is set at $30$ particles per meter.

Simulations are run with both the naive (power) and explog forms for the acoustic term,
and we look at the power density spectra for the pressure oscillations at the prescribed probes.
As shown in Figure~\ref{fig:fft} (bottom) for Probe~6,
the explog choice leads to lower frequency oscillations
($10$ to $20$dB/Hz less throughout most of the spectrum),
and is missing the large spikes that can be observed in the pressure plot
for the naive formulation (Figure~\ref{fig:fft}, top).

\paragraph{Standing wave}

\begin{figure}
\includegraphics[width=\columnwidth]{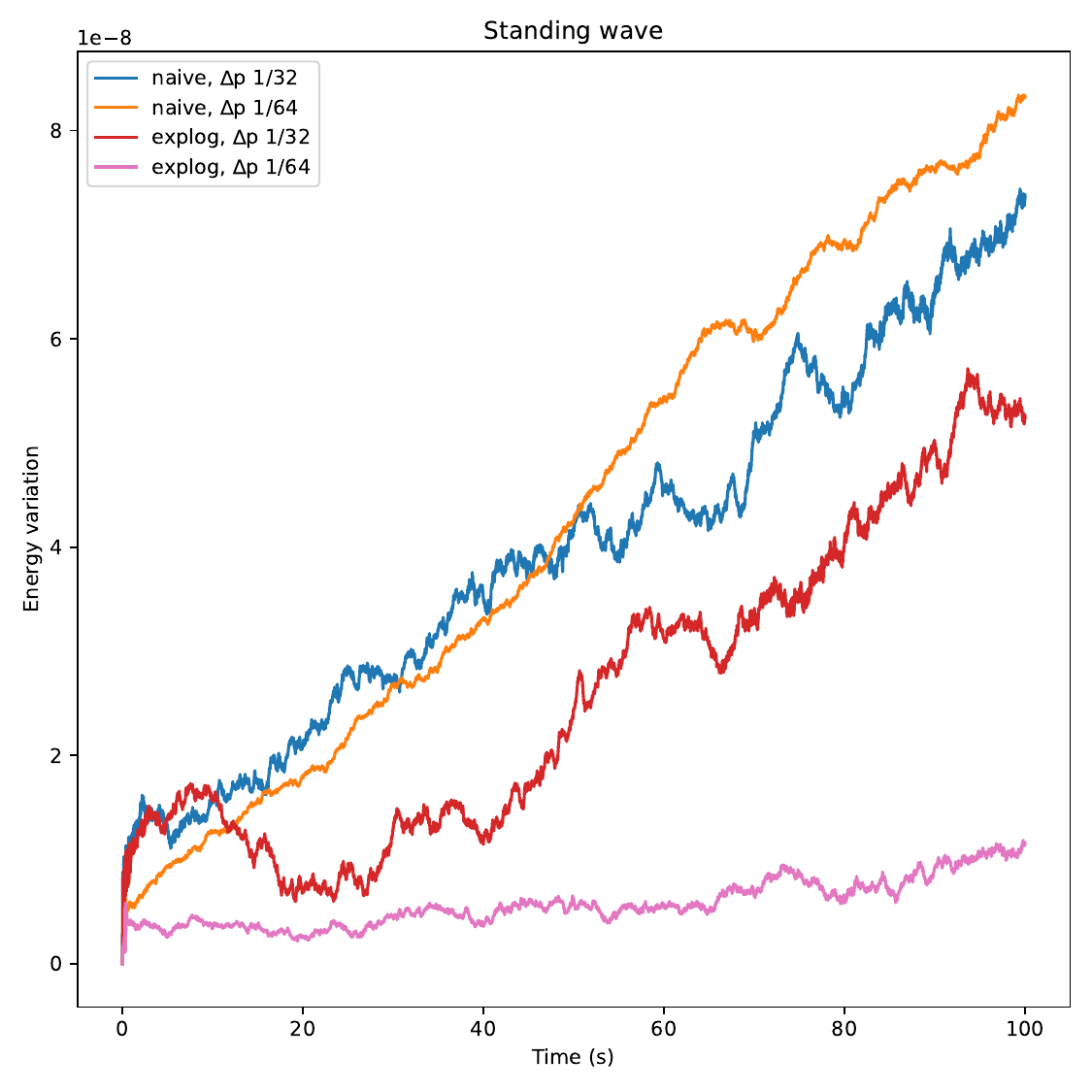}
\caption{Energy conservation in a long-running Standing Wave,
using the naive and the explog forms of the acoustic terms,
and two different particle spacings of $\Delta p = 1/32$
and $\Delta p = 1/64$ particles per meter.}
\label{fig:standing}
\end{figure}

In this test case, we look at energy conservation for a long-running Standing Wave simulation.
The simulation is run in 3D, but with a pseudo-2D configuration.
The domain box is $1$m wide in the wave direction $x$, and $1.1$m tall;
in the cross-wave direction $y$, the domain is only long enough to cover 3 influence radii.
Periodic boundary conditions are applied in both the $x$ and $y$ directions.
Solid boundaries are implemented using the dynamic boundary model,
and the conservative corrected SPH formulation introduced in~\cite{zago_ccsph} is employed
to avoid the well-known numerical dissipative effect in standard WCSPH.
The initial water level is constant, but particles are initialized with non-zero velocity components
$u_x, u_z$ to realize a wave with wavelength $L = 1$m and amplitude $A = 1/64$m.

We run simulations at both a low resolution of 32 particles per meter
and a high resolution of 64 particles per meter,
again comparing the naive (power) and explog forms for the acoustic terms.

We observe (Figure~\ref{fig:standing}) that the naive expression
shows a small but consistent linearly growing trend for the total energy,
and that increasing the resolution makes the energy conservation \emph{worse}
rather than better.
Note that, while the total energy variation remains below $10^{-7}$ in this case,
the trend would push it to grow over single-precision machine epsilon
in longer simulations, with the issue manifesting \emph{earlier} at finer resolutions.

When using the explog form, energy conservation is better overall,
and it improves at finer resolutions, as expected for a numerically robust algorithm.

\paragraph{Stillwater with Irregular topography}

The stillwater case is a well-known benchmark used to evaluate the capability of different density diffusion options and boundary condition formulations to reproduce the pressure field, since its results can be directly compared with the analytical solution.
Prior studies have examined the pressure distribution in a small box of stillwater over a flat bed~\cite{adami_2012, antuono2012stillwater,omidvar2012wavebody,antuono2010freeSurface} and a bed with a wedge (sharp corner)~\cite{adami_2012, colagrossi2012particlepacking, english2022mDBC, fourtakas2019stillwater},
showing that the dynamic boundary condition~\cite{Crespo:2007, dalrymple2001sph} fails to accurately capture the pressure near solid boundaries,
and that modifications to the boundary treatment and density diffusion approach are necessary to improve the results.
Moreover, Colagrossi \emph{et al.}~\cite{colagrossi2012particlepacking} investigated the influence of initial particle positioning on the acoustic noise generated during the early stages of the simulation,
and introduced the Particle Packing Algorithm, designed to generate high-quality initial particle distributions to reduce spurious currents.
 
\begin{figure}
\includegraphics[width=\columnwidth]{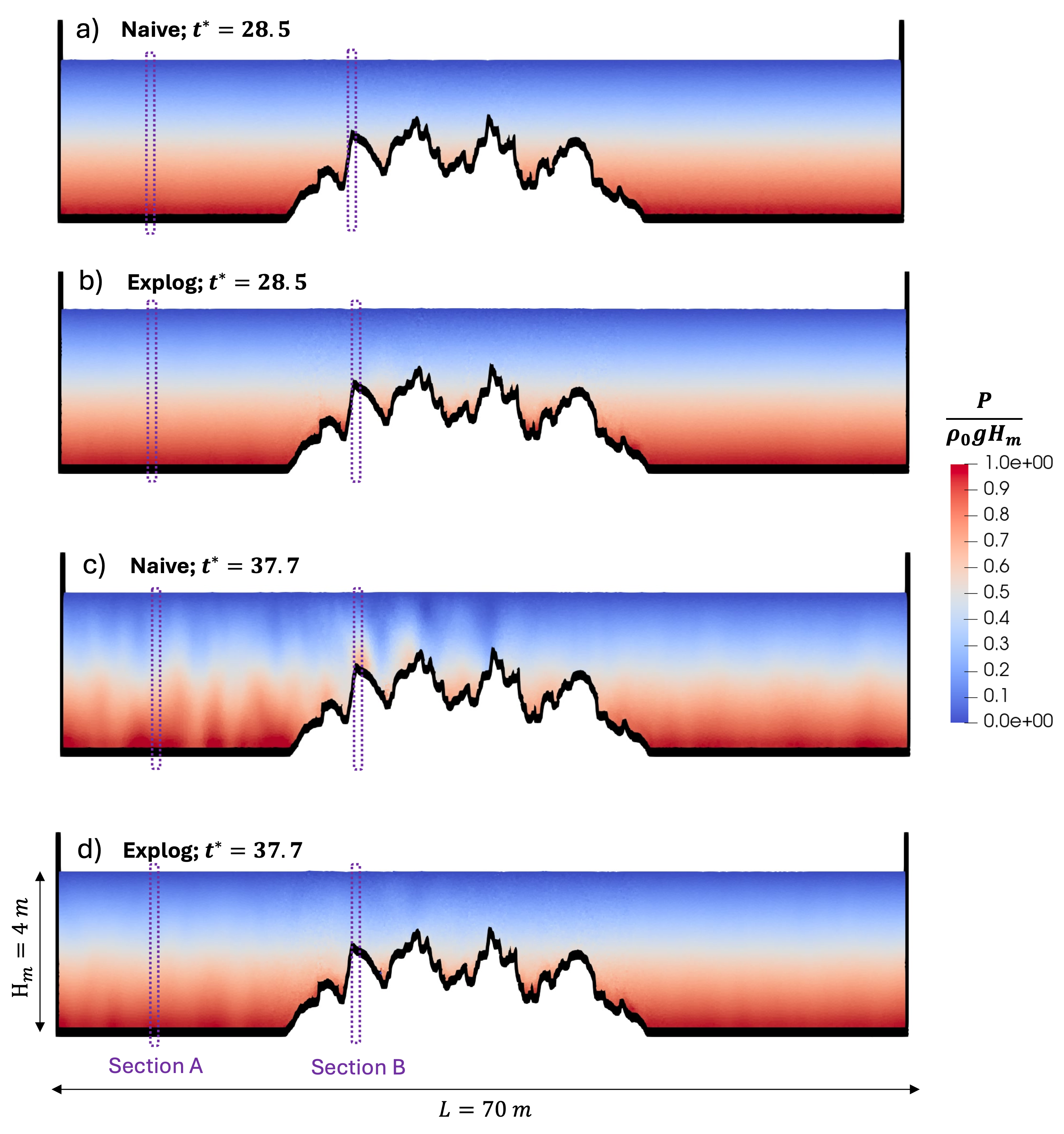}
\caption{Instantaneous normalized pressure fields, \(P/(\rho_0 g H_m)\), for the still-water test: (a) naive and (b) explog at \(t^{*}=28.5\); (c) naive and (d) explog at \(t^{*}=37.7\), where \(t^{*}=t\sqrt{g/H_m}\). All panels use the same bathymetry and color range, \(0\leq (P/\rho_0 g H_m)\leq1\). Purple dotted boxes mark the \(1\)-m-wide particle-extraction windows for Section A (\(15\leq x\leq16~\mathrm{m}\)) and Section B (\(25\leq x\leq26~\mathrm{m}\)), which are used to extract the vertical pressure profiles in Figure~\ref{fig:ptime}.}
\label{fig:Pt}
\end{figure}

To evaluate the capability of the explog formulation in computing pressure distribution near boundaries, particularly in regions with strong bottom roughness,
we model a stillwater testcase where the bottom topography uses a high-resolution bathymetric dataset with rapidly changing bed elevation along the $x$~direction.
As shown in Figure~\ref{fig:Pt}, the irregular DEM consists of two flat sections and a rough segment characterized by steep and rapidly changing topography, including a maximum local slope of \(89^\circ\), and a peak elevation of $2.6$m.
This DEM was used to simulate a still-water case with a maximum water depth of 4.0 m ($H_m = 4\,\text{m}$),
enabling the examination of pressure fluctuations near the bed in the absence of waves and their comparison with the analytical hydrostatic solution over a total simulation time of $100$s.
This duration is substantially longer than that used in previous studies, which typically employed simulation times of less than $10$s~\cite{adami_2012, english2022mDBC, fourtakas2019stillwater}.

Simulations were conducted at a relatively low resolution of 20 particles per meter. Dummy boundary condition was applied to the solid bottom and periodic boundary condition was applied in the alongshore $y$ direction. In order to avoid the excessive nonphysical energy dissipation in the model,
the CCSPH option introduced in~\cite{zago_ccsph} was enabled.
Additionally, the $\delta$-SPH model described in~\cite{marrone_2011} was used for the density diffusion
to reduce acoustic noise and smooth the pressure field.
As shown in Figure \ref{fig:Pt}, pressure distribution was analyzed over time and along the water column at two locations:
one in the flat-bed region (Section A: $x=15$m to $x=16$m),
and another adjacent to the rough-bottom region (Section B: $x=25$m to $x=26$m).
Two different formulations were tested:
the naive formulation where the dummy wall boundary pressure is computed using \eqref{eq:dummy_pressure} directly,
the pressure contribution to the Navier--Stokes equation is computed as $P/\rho^2$ without simplifications,
and the equation of state is computed naively using \code{\_\_powf},
and the explog formulation, with the pressure contribution computed using equations~\eqref{eq:dummy_tilde}--\eqref{eq:dummy_Prhok_tilde}
and the explog form used for $\powipmi$.

Figure~\ref{fig:Pt} compares the normalized pressure fields, \(P/(\rho_0 g H_m)\),
obtained with the naive and explog formulations at two normalized time instances,
\(t^{*}=t\sqrt{g/H_m}\). At \(t^{*}=28.5\) (Figure~\ref{fig:Pt}a,b), both solutions appear nearly steady.
In the naive formulation, however, fluid--boundary interactions over the rough topography generate a localized pressure pulse as numerical errors accumulate.
By \(t^{*}=37.7\) (Figure~\ref{fig:Pt}c)
this pulse has evolved into noticeable pressure fluctuations that originate in the rough section and progressively propagate into the adjacent flat regions.
At the same time, the explog formulation (Figure~\ref{fig:Pt}d) remains nearly steady and does not exhibit comparable fluctuations,
indicating that it limits the growth of the density-related numerical instability triggered by fluid--boundary interactions over the rough topography.

\begin{figure}
\includegraphics[width=0.9\columnwidth]{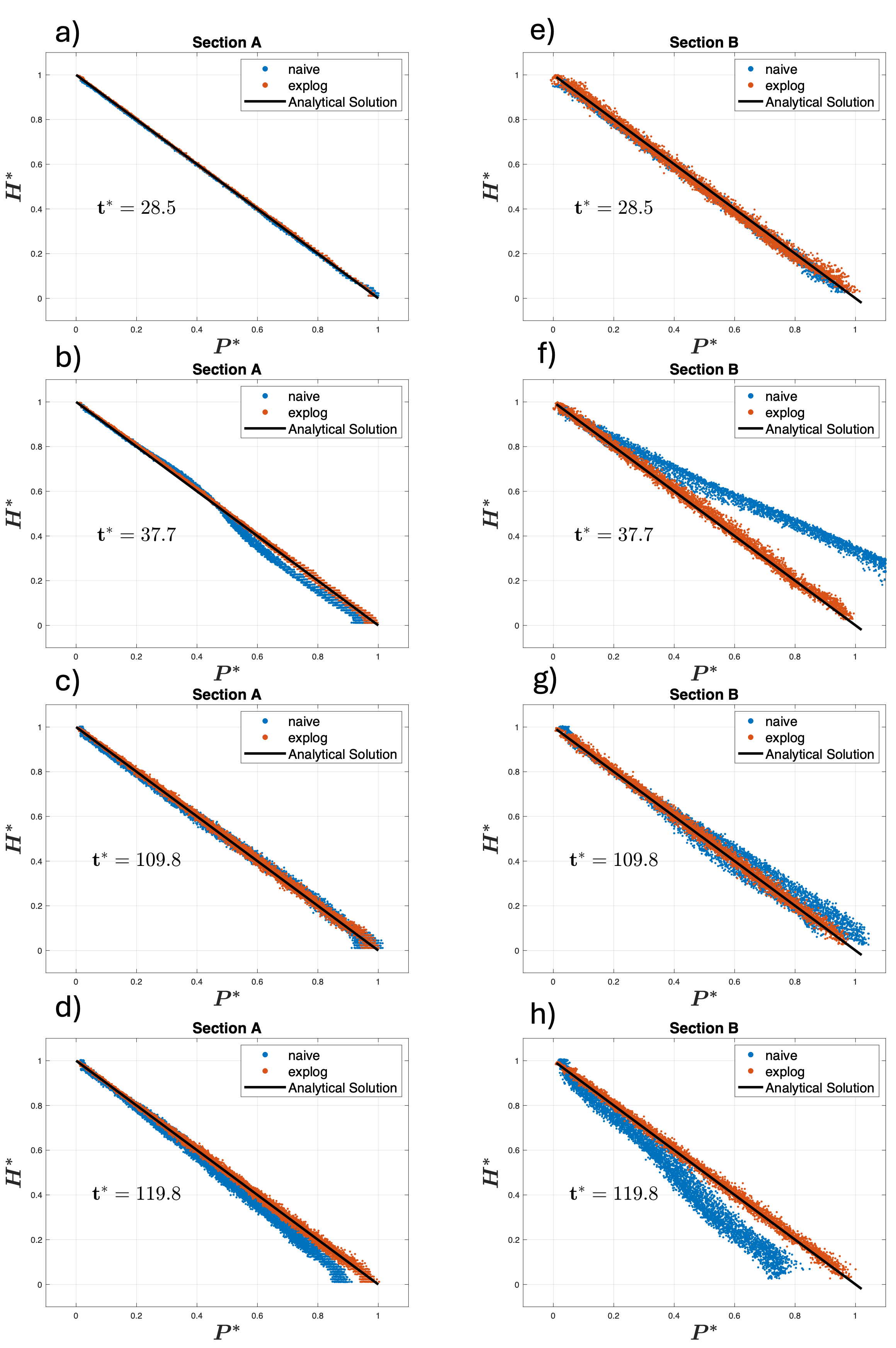}
\caption{Normalized instantaneous particle-pressure profiles in Section A [panels (a)--(d)] and Section B [panels (e)--(h)]. Panels (a,e), (b,f), (c,g), and (d,h) correspond to \(t^{*}=28.5\), \(37.7\), \(109.8\), and \(119.8\), respectively. Blue and orange markers denote fluid-particle pressures computed using the naive and explog formulations, respectively, and the solid black line denotes the analytical hydrostatic solution. Each marker represents one fluid particle within the corresponding \(1\)-m-wide extraction window.}
\label{fig:ptime}
\end{figure}

To quantify these pressure fluctuations,
instantaneous fluid-particle pressures extracted from the \(1\)-m-wide windows in Sections A and B are shown as vertical scatter profiles in Figure~\ref{fig:ptime}. The vertical coordinate \(z\) is measured positive upward from the flat-bed datum,
such that the flat bed in Section A is at \(z=0\) and the undisturbed free surface is at \(z=H_m=4.0~\mathrm{m}\).
For each profile, pressure and vertical position are normalized as \(P_i^{*}=P_i/(\rho_0 g H_{\mathrm{ref}})\)
and \(H_i^{*}=(z_i-z_{\mathrm{b}})/H_{\mathrm{ref}}\), where \(z_{\mathrm{b}}\) is the section reference bottom elevation.
In Section A, \(z_{\mathrm{b}}=0\) and \(H_{\mathrm{ref}}=H_m-z_{\mathrm{b}}=4.0~\mathrm{m}\),
whereas in Section B, \(z_{\mathrm{b}}=2.02~\mathrm{m}\) is the minimum bed elevation within the extraction window,
and \(H_{\mathrm{ref}}=H_m-z_{\mathrm{b}}=1.98~\mathrm{m}\) is the corresponding maximum local still-water depth.
Using the linear hydrostatic reference, the normalized residual for particle \(i\) is defined as
\begin{equation}
e_i^{*}=\frac{P_i-P_{\mathrm{ref}}(z_i)}{\rho_0 g H_{\mathrm{ref}}}, \qquad P_{\mathrm{ref}}(z_i)=\rho_0 g(H_m-z_i),
\label{eq:pressure_residual}
\end{equation}
where \(P_i\) denotes the instantaneous pressure computed by GPUSPH for particle \(i\),
and \(P_{\mathrm{ref}}(z_i)\) is the analytical hydrostatic pressure at the location of particle \(i\).
For each pressure profile, the total departure from the analytical hydrostatic solution and the variability of the pressure residuals are quantified using the normalized root-mean-square error,
\(\mathrm{RMSE}^{*}\), and the normalized residual standard deviation, \(\mathrm{SD}^{*}\), respectively:
\begin{equation}
\begin{aligned}
\mathrm{RMSE}^{*} &= \sqrt{\frac{1}{N}\sum_{i=1}^{N}(e_i^{*})^2}, \\
\mathrm{SD}^{*} &= \sqrt{\frac{1}{N-1}\sum_{i=1}^{N}\left(e_i^{*}-\bar{e}^{*}\right)^2},
\end{aligned}
\label{eq:pressure_metrics}
\end{equation}
where \(N\) is the number of extracted fluid particles and \(\bar{e}^{*}=\frac{1}{N}\sum_{i=1}^{N}e_i^{*}\) is the mean normalized residual.
The \(\mathrm{RMSE}^{*}\) quantifies the overall mismatch between the particle pressures and the analytical hydrostatic solution,
whereas \(\mathrm{SD}^{*}\) quantifies the spread of the residuals about their mean and is used here as a bulk measure of pressure fluctuations.
Reporting both metrics distinguishes a smoothly biased profile, which may have a large \(\mathrm{RMSE}^{*}\) but a small \(\mathrm{SD}^{*}\),
from a profile with substantial residual scatter.

\begin{table}
\centering
\caption{Normalized pressure-residual metrics for the vertical profiles in Sections A and B. \(\mathrm{RMSE}^{*}\) quantifies the total departure of individual particle pressures from the analytical hydrostatic solution, whereas \(\mathrm{SD}^{*}\) quantifies the spread of the normalized pressure residuals about their mean. Values are percentages.}
\label{tab:pressure_metrics}
\small
\setlength{\tabcolsep}{5pt}
\renewcommand{\arraystretch}{1.12}

\begin{tabular}{cccrrrr}
\toprule
& & &
\multicolumn{2}{c}{\(\mathrm{RMSE}^{*}\) (\%)} &
\multicolumn{2}{c}{\(\mathrm{SD}^{*}\) (\%)} \\
\cmidrule(r){4-5}\cmidrule(l){6-7}
\(t\) (s) & \(t^{*}\) & Section & Naive & Explog & Naive & Explog \\
\midrule
\multirow{2}{*}{18.2} & \multirow{2}{*}{28.5} & A & 0.55 & 0.44 & 0.49 & 0.34 \\
& & B & 1.64 & 1.94 & 1.28 & 1.93 \\
\addlinespace[2pt]
\multirow{2}{*}{24.1} & \multirow{2}{*}{37.7} & A & 3.49 & 1.06 & 3.01 & 0.83 \\
& & B & 26.65 & 1.98 & 16.11 & 1.75 \\
\addlinespace[2pt]
\multirow{2}{*}{70.1} & \multirow{2}{*}{109.8} & A & 1.67 & 1.12 & 1.67 & 1.04 \\
& & B & 4.20 & 1.64 & 3.74 & 1.54 \\
\addlinespace[2pt]
\multirow{2}{*}{76.5} & \multirow{2}{*}{119.8} & A & 4.63 & 1.05 & 3.09 & 1.04 \\
& & B & 11.17 & 1.59 & 6.25 & 1.43 \\
\bottomrule
\end{tabular}
\end{table}

Figure~\ref{fig:ptime} and Table~\ref{tab:pressure_metrics} show that the naive formulation does not remain hydrostatic over the full simulation.
At \(t^{*}=28.5\), both formulations are close to the analytical hydrostatic solution with both \(\mathrm{RMSE}^{*}\) and \(\mathrm{SD}^{*}\) below \(2.0\%\).
Figures~\ref{fig:ptime}a and \ref{fig:ptime}e show that, consistent with previous studies,
the model agrees well with the analytical hydrostatic solution over short simulation times.
However, at the next sampled time, the naive formulation begins to deviate from the correct solution,
and this behavior repeats multiple times throughout the simulation.
At \(t^{*}=37.7\), the naive solution departs substantially from the analytical hydrostatic solution,
particularly in Section B [Figure~\ref{fig:ptime}(f)],
where \(\mathrm{RMSE}^{*}\) decreases from \(26.65\%\) for the naive formulation to \(1.98\%\) for explog,
while \(\mathrm{SD}^{*}\) decreases from \(16.11\%\) to \(1.75\%\).
The naive profiles move closer to the analytical solution at \(t^{*}=109.8\),
but the instability reappears at \(t^{*}=119.8\);
in Section B, the naive \(\mathrm{RMSE}^{*}\) and \(\mathrm{SD}^{*}\) are \(11.17\%\) and \(6.25\%\),
respectively, compared with \(1.59\%\) and \(1.43\%\) for explog.
Across all selected times and sections, the explog formulation remains below \(2.0\%\) in \(\mathrm{RMSE}^{*}\) and \(\mathrm{SD}^{*}\). 

Overall, the naive formulation does not consistently reproduce the analytical hydrostatic solution over long simulations,
with the largest departures occurring in Section B near the steep rough topography.
The explog formulation suppresses the growth of the instability observed with the naive formulation,
and remains below \(2.0\%\) in both \(\mathrm{RMSE}^{*}\) and \(\mathrm{SD}^{*}\) at all selected times and sections,
although small residual pressure fluctuations remain.

\subsection{Performance}

We leverage the recently introduced support for CPUs and AMD GPUs in GPUSPH~\cite{bilotta_cpusph}
to test the impact of the choice of expression for $\powipmi$ across a variety of
hardware platforms.

The hardware used for these tests, done again on the Standing Wave test case, include:
\begin{description}
\item[CUDA] an NVIDIA RTX 3060 GPU (Compute Capability 8.6) using CUDA 11.6;
\item[HIP] an AMD Radeon RX 580 GPU (Graphics IP 8) using HIP/ROCm 5.4;
\item[AMD] the AMD Ryzen 7 3700X CPU used for the microbenchmarks presented in
Tables~\ref{tab:pow1pm1.log.1}--\ref{tab:pow1pm1.rt.log.5};
\item[Intel] a dual-socket Intel Xeon Silver 4110 CPU (32 threads total);
\item[ARM] an 8-core Cortex-A53 CPU.
\end{description}

Since we are interested in the possible performance loss (or gain)
related to the use of the explog form,
we compute the absolute performance in MIPPS
(millions of iterations times particles per second, as detailed in~\cite{rustico_advances_2014})
for each configuration, taking the median across~5 tests
(in benchmarking mode, with data output completely disabled to avoid pollution of the results due to I/O latency),
and then express our results as the ratio of the performance
for a given configuration to the performance of the naive (power form) case on the same hardware.
Less than $1\%$ variation was observed in the timings across runs for each configuration.

\begin{figure*}
\centerline{%
\includegraphics[width=0.5\columnwidth]{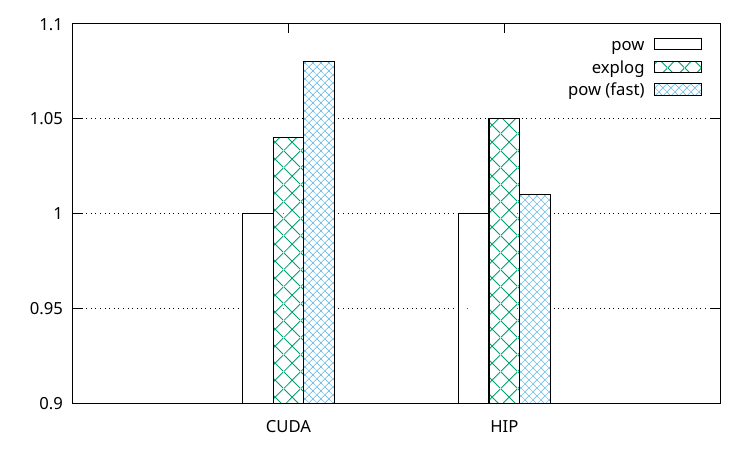}\hfill
\includegraphics[width=0.5\columnwidth]{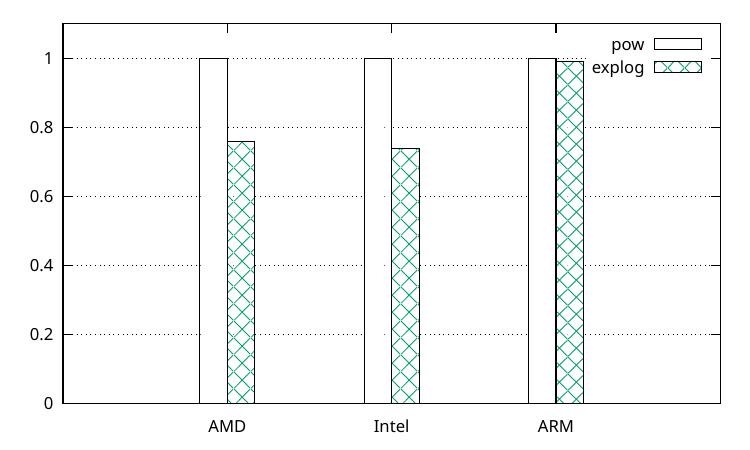}%
}
\caption{Relative performance of the explog and ``fast power'' expressions for $\powipmi$
across different GPUs (left) and CPUs (right). Note that the GPU plot uses a truncated scale to highlight
the relatively small performance difference.}
\label{fig:standing-perf}
\end{figure*}

The results are shown in Figure~\ref{fig:standing-perf}, separately for GPUs and for CPUs.
Since both CUDA and HIP define a lower-accuracy, faster version of the power function (\code{\_\_powf}),
the timing results leveraging this form are also included for these platforms.

On GPUs, we observe that the explog form is not only more accurate, but also marginally
faster ($4\%$ to $5\%$, depending on platform).
The faster power variant achieved even higher performance on CUDA (as expected),
but not on HIP, probably due to differences in the way this is implemented on the AMD platform.

The small performance gains observed on GPU are consistent with the microbenchmarks,
that show the explog form to be between $30\%$ faster
(Table~\ref{tab:pow1pm1.log.1}) to less than $20\%$ slower (Table~\ref{tab:pow1pm1.log.5}),
depending on argument magnitude,
when considering that in a typical WCSPH simulations the arguments will typically be closer to~$0$.

The CPU results are more surprising in this regard:
although the choice of formulation is mostly indifferent on the ARMv8 architecture,
the explog performance is consistently worse on the x86-64 architectures,
with a drop between $20\%$ and $25\%$.
The reasons for these discrepancies from the microbenchmarks are still under analysis,
with higher register consumption being a probable candidate.

Given the hardware-dependent performance ratio between the \code{explog} and naive formulations,
in GPUSPH we have left the choice to the user, who can then balance their computational
and accuracy requirements based on the specific test cases that they need to run.

\section{Conclusions}

Adopting the relative density variation (condensation) and a scaled pressure as state variables
leads to SPH formulations that are numerically more stable.
Additionally, it allows expressing the acoustic terms derived from Cole's equation of state
as a combination of the \code{expm1} and \code{log1p} IEEE-754 functions;
the resulting explog form achieves higher precision for (relative) values closer to zero,
which is the range where condensation and scaled pressure are mostly distributed
in a standard WCSPH simulation.
The explog form of the acoustic terms can also be leveraged to improve pressure evaluation
for boundary terms, which helps reduce spurious pressure oscillations,
and improves energy conservation.

While analytical strategies such as artificial viscosity and density diffusion can help diffuse numerical noise,
they are insufficient in preventing the accumulation of noise over time,
leading to degrading energy conservation and loss of accuracy in long-running simulations.
Combining them with the explog form we introduced ensures that
energy is conserved to below the machine epsilon limit
and hydrostatic density profiles are preserved even for simulations running for hundreds of seconds.

Performance-wise,
the new form also provides performance benefits over the IEEE-compliant power function form on GPU,
but is slower than the non-compliant and less accurate \code{\_\_powf} on CUDA.
On CPU, our implementation of the explog form results in a measurable loss in performance,
that is currently being investigated.

\section*{CRediT}
\textbf{Giuseppe Bilotta:}
Conceptualization, Methodology, Software, Validation, Formal Analysis,
Investigation, Writing - Original Draft, Writing - Review \&\ Editing, Supervision, Project administration, Funding acqusition
\textbf{Alireza Zarei:}
Software, Validation, Investigation, Writing - Review \&\ Editing, Visualization
\textbf{Roberta Cristofaro:}
Investigation, Writing - Review \&\ Editing, Visualization
\textbf{Alexis Hérault:}
Methodology, Software, Formal Analysis, Writing - Review \&\ Editing
\textbf{Élie Saikali:}
Methodology, Software, Formal Analysis, Writing - Review \&\ Editing
\textbf{Morteza Derakhti:}
Validation, Writing - Review \&\ Editing, Supervision, Project administration

\section*{Acknowledgements}
This work was supported by the INGV HARVEST project (BRL~2021, Ob.Fu.~9999.804);
and the 2022 Progetto di Ricerca di Rilevante Interesse Nazionale (PRIN~2022)
\emph{PLastic Transport due to waves and currents ON Emerged and submerged beaches}
(CUP~D53D23004590006, Ob.Fu.~1216.010).

\bibliographystyle{elsarticle-num}
\bibliography{gpusph,nongpusph}

\end{document}